\documentclass[11pt,a4paper,reqno,english]{amsart}

\usepackage{latexsym}
\usepackage{graphicx}
\usepackage{mathrsfs}
\usepackage{amsfonts}
\usepackage{amssymb}
\usepackage{amsmath}
\usepackage{amsthm}
\usepackage{amscd}
\usepackage{color}
\usepackage{xcolor}
\usepackage[all,2cell]{xy}
\usepackage[pagebackref,colorlinks]{hyperref}
\usepackage{enumerate}
\usepackage{bm}
\usepackage{bbm}
\usepackage[normalem]{ulem}
\usepackage{mathrsfs}
\usepackage{rotating}
\usepackage{comment}
\usepackage{hyperref}
\usepackage{comment}
\usepackage{mathtools}
\usepackage{tikz}
\makeatletter
\newcommand\mathcircled[1]{%
  \mathpalette\@mathcircled{#1}%
}
\newcommand\@mathcircled[2]{%
  \tikz[baseline=(math.base)] \node[draw,circle,inner sep=1pt] (math) {$\m@th#1#2$};%
}
\makeatother

\numberwithin{equation}{section}

\theoremstyle{plain}
\newtheorem{thm}{Theorem}[section]
\newtheorem{prop}[thm]{Proposition}
\newtheorem{cor}[thm]{Corollary}
\newtheorem{lemma}[thm]{Lemma}

\theoremstyle{definition}
\newtheorem{deff}[thm]{Definition}
\newtheorem{example}[thm]{Example}

\theoremstyle{remark}
\newtheorem{rmk}[thm]{\bf Remark}

\def\N{\mathbb N}
\def\l{\langle}
\def\r{\rangle}

\newcommand{\ima}{\operatorname{im}}

\newcommand{\reg}{\operatorname{reg}}

\newcommand{\Mat}{\operatorname{\mathbb{M}}}

\newcommand{\V}{\mathcal{V}}
\newcommand{\diag}{\operatorname{diag}}
\newcommand{\supp}{\operatorname{supp}}

\newcommand{\id}{\operatorname{id}}
\newcommand{\blue}{\operatorname{blue}}
\newcommand{\red}{\operatorname{red}}
\newcommand{\base}{\operatorname{base}}
\newcommand{\mor}{\operatorname{\sim_{Mor}}}
\newcommand{\Mid}{\operatorname{\hspace{0.07cm}\big|\hspace{0.11cm}}}

\begin{document}

\title[Moves for Bergman algebras]{Moves for Bergman algebras}

\author{R. Preusser}
\address{Raimund Preusser: School of Mathematics and Statistics\\
Nanjing University of Information Science and Technology, Nanjing, China}
\email{raimund.preusser@gmx.de}



\begin{abstract}
We define Bergman presentations and Bergman algebras associated to Bergman presentations. These algebras embrace various generalisations of Leavitt path algebras. A Bergman presentation can be visualised by a Bergman graph, which is a finite bicoloured hypergraph satisfying two conditions. We define several moves for Bergman graphs and prove that they preserve the isomorphism class (respectively the Morita equivalence class) of the corresponding Bergman algebra. One recovers the well-known results, that in the context of finite directed graphs the shift move, outsplitting, insplitting, source elimination and collapsing preserve the isomorphism class (respectively the Morita equivalence class) of the corresponding Leavitt path algebra. Moreover, we mention some connections between Tietze transformations and the moves for Bergman graphs defined in this paper.
\end{abstract}

\maketitle


\section{Introduction}
Leavitt path algebras of (directed) graphs were introduced by G. Abrams and G. Aranda Pino in 2005 \cite{AA} and independently by P. Ara, M. Moreno and E. Pardo in 2007 \cite{AMP}. The Leavitt path algebras turned out to be a very rich and interesting class of algebras, whose studies so far have comprised over 200 research papers. The field of Leavitt path algebras is a very active research area with connections to functional analysis, symbolic dynamics, K-theory and noncommutative geometry. A comprehensive treatment of the subject can be found in the book~\cite{lpabook}. 

The definition of the Leavitt path algebras was inspired by the algebras $L(m,n)$ studied by W. Leavitt in the 1950's and 60's \cite{vitt56, vitt57, vitt62, vitt65}. For positive integers $m<n$, the Leavitt algebra $L(m,n)$ is universal with the property that $L(m,n)^m\cong L(m,n)^{n}$ as left $L(m,n)$-modules. 
The Leavitt path algebras embrace the algebras $L(1,n)$, but not the algebras $L(m,n)$ where $m>1$. There have been several attempts to introduce a generalisation of the Leavitt path algebras that covers all of the algebras $L(m,n)$. In 2012, Ara and Goodearl introduced Leavitt path algebras of {\it separated graphs}, generalising the Leavitt path algebras of graphs \cite{aragoodearl}. One recovers the algebras $L(m,n)$ as corner rings of Leavitt path algebras of separated graphs. In 2013, R. Hazrat introduced Leavitt path algebras of {\it weighted graphs}, which simultaneously generalise the Leavitt path algebras of graphs and the algebras $L(m,n)$ ~\cite{H-1}. In 2020, the author of this paper introduced Leavitt path algebras of {\it hypergraphs}, which cover the Leavitt path algebras of separated graphs and a subclass of the Leavitt path algebras of weighted graphs (namely the Leavitt path algebras of {\it vertex-weighted graphs}) \cite{Raimund2}. Finally in 2021, R. Mohan and B. Suhas introduced Leavitt path algebras of {\it bi-separated graphs}, which embrace all classes of algebras mentioned earlier in these paragraph \cite{mohan-suhas}.

The questions Q1 and Q2 below belong to the most important question for Leavitt path algebras.
\begin{align*}
{\bf Q1:}&\,\,\text{When do two graphs define isomorphic Leavitt path algebras?}\\
\medskip
{\bf Q2:}&\,\,\text{When do two graphs define Morita equivalent Leavitt path algebras?}
\end{align*}
In this generality these questions are still open (if one considers only finite and acyclic graphs, then the answers to Q1 and Q2 are known). But some ``moves'' for graphs are known to preserve the isomorphism class or the Morita equivalence class of the corresponding Leavitt path algebra. For example, the shift move and outsplitting preserve the isomorphism class of the corresponding Leavitt path algebra, while insplitting, source elimination and collapsing preserve the Morita equivalence class of the corresponding Leavitt path algebra (see \cite[Theorem 1.1]{kbgm}, \cite[\S 6.2, \S 6.3]{lpabook} and \cite[\S 2]{AN}). To the best of our knowledge, for separated graphs, weighted graphs, hypergraphs or bi-separated graphs such moves (preserving the isomorphism class or the Morita equivalence class of the corresponding Leavitt path algebra) have not yet been described in the literature.

In this paper we define Bergman presentations. A Bergman presentation is a finite presentation of an abelian monoid together with a bipartition of the relations such that a few conditions are satisfied. We use some of the universal ring constructions introduced by G. Bergman in \cite{bergman74} to define the Bergman algebra $B(X,R)$ of a Bergman presentation $(X,R)$. The algebra $B(X,R)$ has the property that its $\V$-monoid $\V(B(X,R))$ is presented by $(X,R)$ (for the definition of the $\V$-monoid of a ring see \S 2.2). The Bergman algebras $B(X,R)$ embrace the Leavitt path algebras of finite graphs, finite separated graphs, finite weighted graphs, finite hypergraphs and finite bi-seperated graphs, as well as the Cohn-Leavitt path algebras based on finite graphs or finite separated graphs (for more information on Cohn-Leavitt path algebras see \cite{lpabook} and \cite{aragoodearl}). A Bergman presentation can be visualised by a Bergman graph, which is a finite red-blue coloured hypergraph satisfying two conditions. We show that the category ${\bf BP}$ of Bergman presentations is isomorphic to the category ${\bf BG}$ of Bergman graphs. By definition, the Bergman algebra $B(H)$ of a Bergman graph $H$ is the Bergman algebra of the corresponding Bergman presentation. 

We define several ``moves" for Bergman presentations and Bergman graphs, namely the red shift move, blue shift move, enqueuing, outsplitting, lonely generator/vertex elimination, collapsing and insplitting. We prove that the first four of these moves preserve the isomorphism class of the Bergman algebra, and that the last three of them preserve the Morita equivalence class of the Bergman algebra. One recovers the results mentioned in the third paragraph of this introduction, that in the context of usual graphs the shift move and outsplitting preserve the isomorphism class of the corresponding Leavitt path algebra, while source elimination, collapsing and insplitting preserve the Morita equivalence class of the corresponding Leavitt path algebra. Moreover, we show that any collapsing is a composition of a finite number of red shift moves and lonely generator/vertex eliminations, and that any insplitting is a composition of a finite number of red shift moves and extensions, where by an extension we mean the move that is inverse to a lonely generator/vertex elimination. These two results seem to be new even in the context of graphs (a red shift move corresponds to a shift move and a lonely generator/vertex elimination corresponds to a source elimination in this context).  

The rest of this article is organised as follows. In Section 2, we recall some definitions and results which will be used throughout the paper. In Section 3, we define Bergman presentations, Bergman graphs and Bergman algebras. In Section 4, we define and investigate the moves for Bergman presentations and Bergman graphs mentioned above. In Section 5, we explore connections between these moves and Tietze transformations.

\section{Preliminaries}
Throughout the paper $K$ denotes a fixed field. All rings and algebras are associative and unital unless otherwise stated. By a module we mean a left module. $R\mor S$ means that the rings $R$ and $S$ are Morita equivalent.

Let $R$ be a $K$-algebra. An {\it $R$-ring$_K$} is a $K$-algebra $S$ together with a $K$-algebra homomorphism $R\to S$. An {\it $R$-ring$_K$ homomorphism} from an $R$-ring$_K$ $S$ to an $R$-ring$_K$ $T$ is a $K$-algebra homomorphism $S\to T$ such that the obvious diagram commutes. If $M$ is an $R$-module and $S$ is an $R$-ring$_K$, then we often denote the $S$-module $S\otimes_R M$ also by $M$.

$\N$ denotes the set of positive integers and $\N_0$ the set of nonnegative integers.

By a {\it multiset} over a given set $X$ we mean a map $m:X\to \N_0$. For each $x\in X$, $m(x)$ is called the {\it multiplicity} of $x$. The set $\supp(m)=\{x\in X\mid m(x)\neq 0\}$ is called the {\it support} of $m$. If $\supp(m)\neq\emptyset$, we call $m$ {\it nonempty}. If $\supp(m)$ is finite, we may denote $m$ in the form $\{x_1,\dots,x_n\}$ where $x_1,\dots,x_n\in \supp(m)$ and any $x\in \supp(m)$ occurs precisely $m(x)$ times in the sequence $x_1,\dots,x_n$. If $m(x)=1$ for all $x\in \supp(m)$, we call $m$ a {\it set} and identify $m$ with $\supp(m)$. 

\subsection{Presenting algebras using matrices}
Suppose that $I$ and $J$ are sets, that for any $i\in I$, $\sigma^{(i)}$ is a $m_i\times n_i$ matrix whose entries are symbols, and that for any $j\in J$, $\tau^{(j)}$ and $\xi^{(j)}$ are $m'_j\times n'_j$ matrices over the free $K$-algebra generated by the entries of the matrices $\sigma^{(i)}~(i\in I)$. Then by 
\begin{equation*}
\Big\langle\sigma^{(i)}~(i\in I)\mid \tau^{(j)}=\xi^{(j)}~(j\in J)\Big\rangle
\end{equation*}
we mean the algebra
\begin{align*}
\Big\langle\sigma^{(i)}_{pq}~(i\in I,1\leq p \leq m_i,1\leq q \leq n_i)\mid \tau_{pq}^{(j)}=\xi_{pq}^{(j)}~(j\in J,1\leq p \leq m'_j,1\leq q \leq n'_j)\Big\rangle.\label{eq2.1}
\end{align*}

Now suppose that $A$ and $B$ are $K$-algebras, that $I$ and $J$ are sets, that for any $i\in I$, $\sigma^{(i)}$ is a $m_i\times n_i$ matrix whose entries are symbols, and that for any $j\in J$, $\tau^{(j)}$ and $\xi^{(j)}$ are $m'_j\times n'_j$ matrices over the free product of $A$ and the free $K$-algebra generated by the entries of the matrices $\sigma^{(i)}~(i\in I)$. When we write that $B$ can be obtained from $A$ by adjoining the matrices $\sigma^{(i)}~(i\in I)$ and relations $\tau^{(j)}=\xi^{(j)}~(j\in J)$, we mean that $B$ be can be obtained from $A$ by adjoining the generators $\sigma^{(i)}_{pq}~(i\in I,1\leq p \leq m_i,1\leq q \leq n_i)$ and relations $\tau_{pq}^{(j)}=\xi_{pq}^{(j)}~(j\in J,1\leq p \leq m'_j,1\leq q \leq n'_j)$.

\subsection{The $\V$-monoid of a ring}
Recall that the {\it $\V$-monoid} $\V(R)$ of a ring $R$ is the set of all isomorphism classes of finitely generated projective $R$-modules, which becomes an abelian monoid by defining $[P]+[Q]:=[P\oplus Q]$ for any $[P],[Q]\in \V(R)$. 

There is the following alternative description of the monoid $\V(R)$ using matrices (cf. \cite[p.108]{lpabook}). Let $\Mat_\infty(R)$ be the directed union of the rings $\Mat_n(R)~(n\in\N)$, where the transition maps $\Mat_n(R)\rightarrow \Mat_{n+1}(R)$ are given by $\sigma\mapsto \begin{pmatrix}\sigma&0\\0&0\end{pmatrix}$. Let $I(\Mat_\infty(R))$ denote the set of all idempotent elements of $\Mat_\infty(R)$. Define an equivalence relation $\sim$ on $I(\Mat_\infty(R))$ by $\epsilon_1\sim\epsilon_2$ iff there are $\sigma,\sigma'\in \Mat_\infty(R)$ such that $\epsilon_1 = \sigma\sigma'$ and $\epsilon_2 = \sigma'\sigma$. Let $\hat\V(R)$ be the set of all $\sim$-equivalence classes. The set $\hat\V(R)$ becomes an abelian monoid by defining $[\epsilon_1]+[\epsilon_2]:=[\epsilon_1\oplus\epsilon_2]$ for any $[\epsilon_1],[\epsilon_2]\in \hat \V(R)$, where $\epsilon_1\oplus\epsilon_2=\begin{pmatrix}\epsilon_1&0\\0&\epsilon_2\end{pmatrix}$. There  is a monoid isomorphism $\theta:\hat \V(R)\to \V(R)$ which maps an equivalence class $[\epsilon]\in \hat\V(R)$, where $\epsilon\in \Mat_n(R)$, to the isomorphism class $[R^n\epsilon]\in\V(R)$.

If $\epsilon\in I(\Mat_\infty(R))$ and $\theta([\epsilon])=[P]$ where $P$ is some finitely generated projective $R$-module, we say that $\epsilon$ {\it represents $[P]$}. Note that if $\phi:R\to S$ is a ring homomorphism, $P$ a finitely generated projective $R$-module and $[P]$ is represented by $\epsilon\in I(\Mat_\infty(R))$, then $[S\otimes_R P]$ is represented by the matrix $\phi(\epsilon)\in I(\Mat_\infty(S))$ obtained from $\epsilon$ by applying $\phi$ to each entry. Hence the diagram
\[\xymatrix@C=1.3cm@R=1.3cm{\hat \V(R)\ar[r]^{\phi}\ar[d]_{\theta_R}& \hat \V(S)\ar[d]^{\theta_S}\\\V(R)\ar[r]_{S\otimes_R-}& \V(S)}\]
commutes.

We will need Lemma \ref{lemvmon} below in Section 5.
\begin{lemma}\label{lemvmon}
Let $R$ and $S$ be Morita equivalent rings. Then $\V(R)\cong \V(S)$.
\end{lemma}
\begin{proof}
Let $F:R$-${\bf Mod}\to S$-${\bf Mod}$ and $G:S$-${\bf Mod}\to R$-${\bf Mod}$ be mutually inverse category equivalences, where $R$-${\bf Mod}$ denotes the category of $R$-modules and $S$-${\bf Mod}$ denotes the category of $S$-modules. By \cite[Theorem 18.26]{lam} there is an $R$-$S$ bimodule $P$ and an $S$-$R$ bimodule $Q$ such that $F\cong Q\otimes_R -$ and $G\cong  P\otimes _S-$. It follows that the maps $\V(R)\to \V(S),~[M]\mapsto [F(M)]$ and $\V(S)\to \V(R),~[N]\mapsto [G(N)]$ are mutually inverse monoid homomorphisms.
\end{proof}

\subsection{Some universal ring constructions}
In this subsection $R$ denotes a $K$-algebra. Bergman showed that if $P$ is a finitely generated projective $R$-module, then there is a $R$-ring$_K$ $S$ with a universal idempotent endomorphism $e$ of the $S$-module $S\otimes_R P$ \cite[\S 3]{bergman74}. Note that the isomorphism class of $S$ only depends on the isomorphism class of $P$ (see \cite[p.38]{bergman74}), and that $e$ induces a universal direct sum decomposition $S\otimes_R P=P_1\oplus P_2$ where $P_1=\ker(e)$ and $P_2=\ima(e)$. While Bergman denoted $S$ by $R\langle e:\overline{P}\to\overline{P};~e^2=e\rangle$, we denote it by 
\[R\Big\langle [P_1],[P_2]\Mid [P]=[P_1]+[P_2]\Big\rangle.\]

Bergman also showed that if $P$ and $Q$ are finitely generated projective $R$-modules, then there is a $R$-ring$_K$ $T$ with a universal isomorphism $i:T\otimes_RP\to T\otimes_RQ$ of $T$-modules \cite[\S 3]{bergman74}. Note that the isomorphism class of $T$ only depends on the isomorphism classes of $P$ and $Q$ (see \cite[p.38]{bergman74}). While Bergman denoted $T$ by $R\langle i,i^{-1}:\overline{P}\cong\overline{Q}\rangle$, we denote it by \[R\Big\langle [P]=[Q]\Big\rangle.\]

If $t\geq 2$ and $P$ is a finitely generated projective $R$-module, then we denote the $R$-ring$_K$ 
\begin{align*}
R\Big\langle [P_1],[P_2] \Mid [P]=[P_1]+[P_2]\Big\rangle\Big\langle [P'_2],[P_3]\Mid [P_2]=[P'_{2}]+ [P_3]\Big\rangle \dots \Big\langle [P'_{t-1}],[P_t] \Mid [P_{t-1}]=[P'_{t-1}]+ [P_t]\Big\rangle
\end{align*}
by \[R\Big\langle [P_1],\dots,[P_t]\Mid [P]=[P_1]+\dots +[P_t]\Big\rangle.\]

 
\begin{lemma}\label{lembergpres1}
Let $P$ be a finitely generated projective $R$-module and $\zeta\in \Mat_n(R)$ an idempotent matrix representing $[P]$. Then the $K$-algebra $S=R\Big\langle [P_1],[P_2]\Mid [P]=[P_1] +  [P_2]\Big\rangle$ can be obtained from $R$ by adjoining an $n\times n$ matrix $\epsilon$ and the relations $\zeta\epsilon\zeta=\epsilon$ and $\epsilon^2=\epsilon$. If $P$ is nonzero, then $\V(S)$ can be obtained from $\V(R)$ by adjoining two new generators $[P_1]$ and $[P_2]$ and one relation $[P]=[P_1]+[P_2]$.
\end{lemma}
\begin{proof}
The first assertion of the lemma follows from \cite[Proofs of Theorems 3.1 and 3.2]{bergman74}, the second assertion follows from \cite[Theorem 5.1]{bergman74}.
\end{proof}

We leave it to the reader to deduce Lemma \ref{lembergpres1.5} below from Lemma \ref{lembergpres1} above.
\begin{lemma}\label{lembergpres1.5}
Let $t\geq 2$, $P$ a finitely generated projective $R$-module and $\zeta\in \Mat_n(R)$ an idempotent matrix representing $[P]$. Then the $K$-algebra $S=R\Big\langle [P_1],\dots,[P_t]\Mid [P]=[P_1] +  \dots  +  [P_t]\Big\rangle$ can be obtained from $R$ by adjoining $n\times n$ matrices $\epsilon_1,\dots,\epsilon_{t-1}$ and relations $\zeta\epsilon_i\zeta=\epsilon_i$ and $\epsilon_i\epsilon_j=\delta_{ij}\epsilon_i$ where $1\leq i,j\leq t-1$ and $\delta_{ij}$ is the Kronecker delta. If $P$ is nonzero, then $\V(S)$ can be obtained from $\V(R)$ by adjoining $t$ new generators $[P_1],\dots,[P_{t}]$ and one relation $[P]=[P_1]+\dots+[P_{t}]$.
\end{lemma}

\begin{lemma}\label{lembergpres2}
Let $P$ and $Q$ be finitely generated projective $R$-modules, and $\zeta\in \Mat_m(R)$ and $\eta\in \Mat_n(R)$ idempotent matrices representing $P$ and $Q$, respectively. Then the $K$-algebra $T=R\Big\langle [P]=[Q]\Big\rangle$ can be obtained from $R$ by adjoining an $m\times n$ matrix $\sigma$, an $n\times m$ matrix $\sigma'$ and the relations $\zeta\sigma\eta=\sigma$, $\eta\sigma'\zeta=\sigma'$, $\sigma\sigma'=\zeta$ and $\sigma'\sigma=\eta$. If $P$ and $Q$ are nonzero, then $\V(T)$ can be obtained from $\V(R)$ by imposing one relation $[P]=[Q]$.
\end{lemma}
\begin{proof}
The first assertion of the lemma follows from \cite[Proofs of Theorems 3.1 and 3.2]{bergman74}, the second assertion follows from \cite[Theorem 5.2]{bergman74}.
\end{proof}

\begin{lemma}\label{lemoutsplit}
Let $t\geq 2$ and $P_1,\dots,P_t, P,Q$ be finitely generated projective $R$-modules such that $[P]=[P_1] + \dots  +  [P_t]$ in $\V(R)$. Then 
\[R\Big\l [P]=[Q]\Big\r\cong R\Big\l [Q_1],\dots,[Q_t]\Mid [Q]=[Q_1] + \dots  +  [Q_t]\Big\r\Big\l [P_1]=[Q_1]\Big\r\dots\Big\l [P_t]=[Q_t]\Big\r.\]
\end{lemma}
\begin{proof}
Let $[P_i]$ be represented by the matrix $\zeta_i\in \Mat_{m_i}(R)$ for any $1\leq i\leq t$. Then $[P]$ is represented by the matrix $\zeta=\zeta_1\oplus\dots\oplus\zeta_t\in \Mat_m(R)$ where $m=m_1+\dots+m_t$. Let $[Q]$ be represented by the matrix $\eta\in \Mat_n(R)$. Set $A:=R\Big\l [P]=[Q]\Big\r$ and $B:=R\Big\l [Q_1],\dots,[Q_t]\Mid [Q]=[Q_1] + \dots  +  [Q_t]\Big\r\Big\l [P_1]=[Q_1]\Big\r\dots\Big\l [P_t]=[Q_t]\Big\r$. By Lemma \ref{lembergpres2}, $A$ can be obtained from $R$ by adjoining an $m\times n$ matrix $\sigma$, an $n\times m$ matrix $\sigma'$ and the relations $\zeta\sigma\eta=\sigma$, $\eta\sigma'\zeta=\sigma'$, $\sigma\sigma'=\zeta$ and $\sigma'\sigma=\eta$. By Lemmas \ref{lembergpres1.5} and \ref{lembergpres2}, $B$ can be obtained from $R$ by adjoining $n\times n$ matrices $\epsilon_1,\dots,\epsilon_{t-1}$, an $m_i\times n$ matrix $\tau_i$ and an $n\times m_i$ matrix $\tau'_i$ for any $1\leq i\leq t$, and relations $\eta\epsilon_i\eta=\epsilon_i~(1\leq i\leq t-1)$, $\epsilon_i\epsilon_j=\delta_{ij}\epsilon_i~(1\leq i,j\leq t-1)$, $\zeta_i\tau_i\epsilon_i=\tau_i~(1\leq i\leq t)$, $\epsilon_i\tau_i'\zeta_i=\tau_i'~(1\leq i\leq t)$, $\tau_i\tau_i'=\zeta_i~(1\leq i\leq t)$ and $\tau_i'\tau_i=\epsilon_i~(1\leq i\leq t)$. Here $\epsilon_t:=\eta-\epsilon_1-\dots-\epsilon_{t-1}$.

Write the matrices $\sigma$ and $\sigma'$ in block form $\sigma=\begin{pmatrix}\sigma_{1}|\dots|\sigma_{t}\end{pmatrix}^T$ and $\sigma'=\begin{pmatrix}\sigma'_{1}|\dots|\sigma'_{t}\end{pmatrix}$ where for any $1\leq i\leq t$, $\sigma_i$ is an $m_i\times n$ matrix and $\sigma'_i$ is an $n\times m_i$ matrix. We leave to the reader to check that there are $R$-ring$_K$ homomorphisms $\alpha:A\to B$ and $\beta:B\to A$ such that $\alpha(\sigma_{i})=\tau_i$, $\alpha(\sigma'_{i})=\tau'_i$, $\beta(\epsilon_i)=\sigma'_{i}\sigma_{i}$, $\beta(\tau_{i})=\sigma_{i}$ and $\beta(\tau'_{i})=\sigma'_{i}$. Clearly these homomorphisms are inverse to each other. Thus $A\cong B$.
\end{proof}

\subsection{Graphs}
A {\it (directed) graph} is a quadruple $E=(E^0,E^1,s,r)$ where $E^0$ and $E^1$ are sets and $s,r:E^1\rightarrow E^0$ maps. The elements of $E^0$ are called {\it vertices} and the elements of $E^1$ {\it edges}. If $e$ is an edge, then $s(e)$ is called its {\it source} and $r(e)$ its {\it range}. In this article all graphs $E$ are assumed to be {\it finite}, i.e. $E^0$ and $E^1$ are finite sets. In this setup a vertex $v$ is called {\it regular} if $s^{-1}(v)\neq\emptyset$. The subset of $E^0$ consisting of all regular vertices is denoted by $E^0_{\reg}$.

A {\it separated graph} is a pair $(E,C)$ where $E$ is a graph, $C=\bigsqcup_{v\in E^0}C_v$ and $C_v$ is a partition of $s^{-1}(v)$ (into pairwise disjoint nonempty subsets) for every vertex $v$.

A {\it weighted graph} is a pair $(E,w)$ where $E$ is a graph and $w:E^1\rightarrow \N$ is a map. If $e\in E^1$, then $w(e)$ is called the {\it weight} of $e$. If $v\in E^0$, then $w(v):=\max\{w(e)\mid e\in s^{-1}(v)\}$ is called the {\it weight} of $v$ (here we use the convention $\max\emptyset=0$). If $w(e)=w(s(e))$ for any $e\in E^1$, then $(E,w)$ is called a {\it vertex-weighted graph}.

A {\it hypergraph} is a quadruple $H=(H^0,H^1,s,r)$, where $H^0$ and $H^1$ are sets and $s$ and $r$ are maps associating to each $h\in H^1$ a nonempty multiset $s(h)$ respectively $r(h)$ over $H^0$. 
The elements of $H^0$ are called {\it vertices} and the elements of $H^1$ {\it hyperedges}. If $h$ is hyperedge, then any vertex $v\in \supp(s(h))$ is called {\it a source of $h$ (with multiplicity $s(h)(v)$)}. Similarly, any vertex $v\in \supp(r(h))$ is called {\it a range of $h$ (with multiplicity $r(h)(v)$)}. In this article all hypergraphs $H$ are assumed to be {\it finite}, i.e. $H^0$ and $H^1$ are finite sets.

A {\it bi-separated graph} is a triple $\dot{E}=(E,C,D)$ where $E$ is a graph, $C=\bigsqcup_{v\in E^0}C_v$ where $C_v$ is a partition of $s^{-1}(v)$ for any $v\in E^0$, and $D=\bigsqcup_{v\in E^0}D_v$ where $D_v$ is a partition of $r^{-1}(v)$ for any $v\in E^0$, such that $|X\cap Y|\leq 1$ for any $X\in C$ and $Y\in D$. For $X\in C$ we denote by $s(X)$ the common source of the edges in $X$, and for $Y\in D$ we denote by $r(Y)$ the common range of the edges in $Y$. Moreover, for $X\in C$ and $Y\in D$ we define $XY=YX=e$ if $X\cap Y=\{e\}$, respectively $XY=YX=0$ if $X\cap Y=\emptyset$.

Let $\dot{E}=(E,C,D)$ be a bi-separated graph. We say that {\it $\dot{E}$ satisfies Condition (B)} if there is an ordering $X_1,\dots, X_m$ of the elements of $C$ and an ordering $Y_1,\dots,Y_n$ of the elements of $D$ such that the $m\times n$ matrix $A$ whose entry at position $(i,j)$ equals $X_iY_j$ has the block diagonal form $A=\diag(A^{(1)},\dots,A^{(p)})$ where for each $1\leq k\leq p$ the (not necessarly quadratic) matrix $A^{(k)}$ has the upper triangular block form
\[A^{(k)}=\begin{pmatrix}
A^{(k)}_{11}&A^{(k)}_{12}&A^{(k)}_{13}&\dots&A^{(k)}_{1q_k}\\0&A^{(k)}_{22}&A^{(k)}_{23}&\dots&A^{(k)}_{2q_k}\\0&0&A^{(k)}_{33}&\dots&A^{(k)}_{3q_k}\\\vdots&\vdots&\ddots&\ddots&\vdots\\0&0&\dots&0&A^{(k)}_{q_kq_k}
\end{pmatrix}\] 
where all entries of the (not necessarly quadratic) matrices $A^{(k)}_{ij}~(1\leq i,j\leq q_k)$ are nonzero.

\begin{rmk}\label{remhyp}
In \cite{Raimund2}, hypergraphs were defined in a different way, using indexed families instead of multisets. In this article we define hypergraphs using multiset, because that makes it easier to describe some of the moves in Section 4.
\end{rmk} 

\subsection{Leavitt path algebras}

\begin{deff}\label{deflpa}
Let $E$ be a graph. The $K$-algebra $L(E)$ presented as a nonunital $K$-algebra by the generating set $\{v,e,e^*\mid v\in E^0,e\in E^1\}$ and the relations
\begin{enumerate}[(i)]
\item $uv=\delta_{uv}u\quad(u,v\in E^0)$,
\medskip
\item $s(e)e=e=er(e),~r(e)e^*=e^*=e^*s(e)\quad(e\in E^1)$,
\medskip
\item $e^*f= \delta_{ef}r(e)\quad(v\in E_{reg}^0;~e,f\in s^{-1}(v))$,
\medskip
\item $\sum\limits_{e\in s^{-1}(v)}ee^*= v\quad(v\in E_{\reg}^0)$
\end{enumerate}
is called the {\it Leavitt path algebra} of $E$. 
\end{deff}

\begin{rmk}\label{remlpa}
Let $E$ be a graph and $R=K^{E^0}$. For any $v\in E^0$ we identify $v$ with the element of $R$ whose $v$-th component is $1$ and whose other components are $0$. For any $v\in E^0_{\reg}$ we define the finitely generated projective $R$-modules $P_v:=Rv$ and $Q_v:=\bigoplus_{e\in s^{-1}(v)}Rr(e)$. Choose any ordering $v_1,\dots, v_n$ of the regular vertices. Then $L(E)\cong R\Big\langle [P_{v_1}]=[Q_{v_1}]\Big\rangle\dots \Big\langle[P_{v_n}]= [Q_{v_n}]\Big\rangle$, see \cite[\S 3.2]{lpabook}.
\end{rmk}

\begin{deff}\label{defslpa}
Let $(E,C)$ be a separated graph. The $K$-algebra $L(E,C)$ presented as a nonunital $K$-algebra by the generating set $\{v,e,e^*\mid v\in E^0,e\in E^1\}$ and the relations
\begin{enumerate}[(i)]
\item $uv=\delta_{uv}u\quad(u,v\in E^0)$,
\medskip
\item $s(e)e=e=er(e),~r(e)e^*=e^*=e^*s(e)\quad(e\in E^1)$,
\medskip
\item $e^*f= \delta_{ef}r(e)\quad(X\in C;~e,f\in X)$,
\medskip
\item $\sum\limits_{e\in X}ee^*= v\quad(X\in C_v,~v\in E_{\reg}^0)$
\end{enumerate}
is called the {\it Leavitt path algebra} of $(E,C)$. 
\end{deff}

\begin{rmk}\label{remslpa}
\begin{enumerate}[(a)]
\item If $(E,C)$ is a separated graph such that $C_v=\{s^{-1}(v)\}$ for all $v\in E_{\reg}^0$, then $L(E,C)\cong L(E)$. Hence the Leavitt path algebras of separated graphs embrace the Leavitt path algebras of graphs.
\item Let $(E,C)$ be a separated graph and $R=K^{E^0}$. For any $v\in E^0$ we identify $v$ with the element of $R$ whose $v$-th component is $1$ and whose other components are $0$. For any $X\in C$ we define the finitely generated projective $R$-modules $P_X:=Rv$, where $v$ is the common source of the edges in $X$, and $Q_X:=\bigoplus_{e\in X}Rr(e)$. Choose any ordering $X_1,\dots, X_n$ of the elements of $C$. Then $L(E,C)\cong R\Big\langle [P_{X_1}]=[Q_{X_1}]\Big\rangle\dots \Big\langle [P_{X_n}]=[Q_{X_n}]\Big\rangle$, see \cite[\S 4]{aragoodearl}.
\end{enumerate}
\end{rmk}

\begin{deff}\label{defwlpa}
Let $(E,w)$ be a weighted graph. The $K$-algebra $L(E,w)$ presented as a nonunital $K$-algebra by the generating set $\{v,e_i,e_i^*\mid v\in E^0, e\in E^1, 1\leq i\leq w(e)\}$ and the relations
\begin{enumerate}[(i)]
\item $uv=\delta_{uv}u\quad(u,v\in E^0)$,
\medskip
\item $s(e)e_i=e_i=e_ir(e),~r(e)e_i^*=e_i^*=e_i^*s(e)\quad(e\in E^1, 1\leq i\leq w(e))$,
\medskip
\item 
$\sum\limits_{1\leq i\leq w(v)}e_i^*f_i= \delta_{ef}r(e)\quad(v\in E_{reg}^0;~e,f\in s^{-1}(v))$,
\medskip 
\item $\sum\limits_{e\in s^{-1}(v)}e_ie_j^*= \delta_{ij}v\quad(v\in E_{\reg}^0;~1\leq i, j\leq w(v))$
\end{enumerate}
is called the {\it Leavitt path algebra} of $(E,w)$. In relations (iii) and (iv) we set $e_i$ and $e_i^*$ zero whenever $i > w(e)$. 
\end{deff}

\begin{rmk}\label{remwlpa}
\begin{enumerate}[(a)]
\item
If $(E,w)$ is a weighted graph such that $w(e)=1$ for all $e \in E^{1}$, then $L(E,w)\cong L(E)$. Hence the Leavitt path algebras of weighted graphs embrace the Leavitt path algebras of graphs.
\item Let $(E,w)$ be a vertex-weighted graph and $R=K^{E^0}$. For any $v\in E^0$ we identify $v$ with the element of $R$ whose $v$-th component is $1$ and whose other components are $0$. For any $v\in E^0_{\reg}$ we define the finitely generated projective $R$-modules $P_v:=\bigoplus_{1\leq i\leq w(v)}Rv$ and $Q_v:=\bigoplus_{e\in s^{-1}(v)}Rr(e)$. Choose any ordering $v_1,\dots, v_n$ of the regular vertices. Then $L(E,w)\cong R\Big\langle [P_{v_1}]=[Q_{v_1}]\Big\rangle\dots \Big\langle [P_{v_n}]=[Q_{v_n}]\Big\rangle$, see \cite[\S 4]{preusser-1}.
\item One can also use the universal ring constructions mentioned in \S 2.3 to describe the Leavitt path algebras of weighted graphs that are not vertex-weighted. But in this case the description is more complicated, see \cite[\S 4]{preusser-1}.
\end{enumerate}
\end{rmk}

\begin{deff}\label{defhlpa}
Let $H$ be a hypergraph. For any hyperedge $h\in H^1$ set $I_h:=\{(u,k)\mid u\in H^0,1\leq k\leq s(h)(u)\}$ and $J_h:=\{(v,l)\mid v\in H^0,1\leq l\leq r(h)(v)\}$. The $K$-algebra $L(H)$ presented as a nonunital $K$-algebra by the generating set $\{v,h_{ij},h_{ij}^*\mid v\in H^0, h\in H^1, i\in I_h,  j\in J_h\}$ and the relations
\begin{enumerate}[(i)]
\item $uv=\delta_{uv}u\quad(u,v\in H^0)$,
\medskip
\item $uh_{ij}=h_{ij}=h_{ij}v,~vh_{ij}^*=h_{ij}^*=h_{ij}^*u\quad(h\in H^1, ~i=(u,k)\in I_h,  ~j=(v,l)\in J_h)$,
\medskip
\item $\sum\limits_{j\in J_h}h_{ij}h_{i'j}^*= \delta_{ii'}u \quad(h\in H^1; ~i=(u,k),i'=(u',k')\in I_h)$ and
\medskip
\item $\sum\limits_{i\in I_h}h_{ij}^*h_{ij'}= \delta_{jj'}v\quad(h\in H^1;~j=(v,l),j'=(v',l')\in J_h)$
\end{enumerate}
is called the {\it Leavitt path algebra of $H$}. 
\end{deff}

\begin{rmk}\label{remhlpa}
\begin{enumerate}[(a)]
\item Let $H$ be a hypergraph. Let $\widetilde s$ be the map associating to each $h\in H^1$ the family $(\widetilde s(h)_i)_{i\in I_h}$, where $I_h$ is defined as in Definition \ref{defhlpa} and $\widetilde s(h)_i=u$ for any $i=(u,k)\in I_h$. Similarly, let $\widetilde r$ be the map associating to each $h\in H^1$ the family $(\widetilde r(h)_j)_{j\in J_h}$, where $J_h$ is defined as in Definition \ref{defhlpa} and $\widetilde r(h)_i=v$ for any $j=(v,j)\in J_h$. Then $\widetilde H=(H^0, H^1, \widetilde s,\widetilde r)$ is a hypergraph in the sense of \cite{Raimund2}, and moreover $L(H)\cong L(\widetilde H)$.
\item The Leavitt path algebras of hypergraphs embrace the Leavitt path algebras of graphs, separated graphs and vertex-weighted graphs. For details see \cite{Raimund2}.
\item Let $H$ be a hypergraph and $R=K^{H^0}$. For any $v\in H^0$ we identify $v$ with the element of $R$ whose $v$-th component is $1$ and whose other components are $0$. For any $h\in H^1$ we define the finitely generated projective $R$-modules $P_h:=\bigoplus_{v\in\supp(s(h))}(\bigoplus_{i=1}^{s(h)(v)}Rv)$ and $Q_h:=\bigoplus_{v\in\supp(r(h))}(\bigoplus_{i=1}^{r(h)(v)}Rv)$. Choose any ordering $h_1,\dots, h_n$ of the hyperedges. Then $L(H)\cong R\Big\langle [P_{h_1}]=[Q_{h_1}]\Big\rangle\dots \Big\langle [P_{h_n}]=[Q_{h_n}]\Big\rangle$, see \cite[\S 9]{Raimund2}.
\end{enumerate}
\end{rmk}

\begin{deff}\label{defbslpa}
Let $\dot E=(E,C,D)$ be a bi-separated graph. The $K$-algebra $L(\dot E)$ presented as a nonunital $K$-algebra by the generating set $\{v,e,e^*\mid v\in E^0,e\in E^1\}$ and the relations
\begin{enumerate}[(i)]
\item $uv=\delta_{uv}u\quad(u,v\in E^0)$,
\medskip
\item $s(e)e=e=er(e),~r(e)e^*=e^*=e^*s(e)\quad(e\in E^1)$,
\medskip
\item $\sum\limits_{Y\in D}(XY)(YX')^*=\delta_{XX'}s(X)\quad (X,X'\in C)$,
\medskip
\item $\sum\limits_{X\in C}(YX)^*(XY')=\delta_{YY'}r(Y)\quad(Y,Y'\in D)$
\end{enumerate}
is called the {\it Leavitt path algebra} of $\dot E$.
\end{deff}

\begin{rmk}\label{rembslpa}

\begin{enumerate}[(a)]
\item
The Leavitt path algebras of bi-separated graphs satisfying Condition (B) embrace the Leavitt path algebras of graphs, separated graphs, weighted graphs and hypergraphs. For details see \cite{mohan-suhas}.
\item One can use the universal ring constructions mentioned in \S 2.3 to describe the Leavitt path algebras of bi-separated graphs satisfying Condition (B), see \cite[\S 6]{mohan-suhas}.
\end{enumerate}
\end{rmk}

\begin{rmk}
Any Leavitt path algebra of a graph, separated graph, weighted graph, hypergraph or bi-separated graph is a unital algebra whose unit is the sum of all vertices (recall that we assume that any graph or hypergraph is finite).
\end{rmk}

\subsection{Linear bases for Leavitt path algebras of hypergraphs}
Let $H$ be a hypergraph and $I_h,J_h~(h\in H^1)$ be the sets defined in Definition \ref{defhlpa}. We define a (usual) graph $E=(E^0,E^1, s',r')$ by $E^0=H^0$, $E^1=\{h_{ij} \mid h\in H^1, i\in I_h,  j\in J_h\}$, $s'(h_{ij})=u$ and $r'(h_{ij})=v$ for any $i=(u,k)\in I_h$ and $j=(v,l)\in J_h$. The graph $E$ is called the {\it graph associated to $H$}. The {\it double graph} $\widehat E$ of $E$ can be obtained from $E$ by adding for any edge $e\in E^1$ an edge $e^*$ with reverse orientation (cf. \cite[Remark 1.2.4]{lpabook}). 

We will need Theorem \ref{thmbasis} below in Section 4. Recall that a {\it path} in a graph $F=(F^0, F^1, s, r)$ is a finite and nonempty word $x_1\dots x_n$ over the alphabet $F^0\cup F^1$ such that either $x_1,\dots,x_n\in F^1$ and $r(f_i)=s(f_{i+1})~(1\leq i\leq n-1)$, or $n=1$ and $x_1\in F^0$.

\begin{thm}\label{thmbasis}
Let $H$ be a hypergraph and $\widehat E$ the double graph of the graph $E$ associated to $H$. For any hyperedge $h\in H^1$ choose an index $i_h\in I_h$ and an index $j_h\in J_h$. The paths in $\widehat E$ that do not contain any of the words 
\[h_{ij_h}h_{i'j_h}^*~(h\in H^1, i,i'\in I_h)\text{ and }h_{i_hj}^*h_{i_hj'}~(h\in H^1,j,j'\in J_h)\]
as a subword form a basis of the $K$-vector space $L(H)$.
\end{thm}
\begin{proof}
The theorem follows from Remark \ref{remhlpa}(a) and \cite[Corollary 19]{Raimund2}.
\end{proof}

\section{Bergman algebras}
In this section we define Bergman presentations, Bergman graphs and Bergman algebras. Moreover, we show that the category ${\bf BP}$ of Bergman presentations is isomorphic to the category ${\bf BG}$ of Bergman graphs.
\subsection{Bergman presentations and Bergman graphs}
\subsubsection{Bergman presentations}
We call a pair $(X,R)$, where $X$ is a set and $R=\{(a_i,b_i)\mid i\in I\}$ is a family of ordered pairs of elements of the free abelian monoid $\l X\r$, an {\it abelian monoid presentation}. We allow $(a_i,b_i)=(a_j,b_j)$ for $i\neq j$. The elements of $X$ are called {\it generators} and the ordered pairs in $R$ {\it relations}. We usually write a relation $(a_i,b_i)$ in the form $a_i=b_i$. The presentation $(X,R)$ is called {\it finite} if $X$ and $I$ are finite, and {\it good} if all $a_i$'s and $b_i$'s are nonzero in $\l X\r$. Note that there is a 1-1 correspondence between finite and good abelian monoid presentations and hypergraphs (cf. \S 3.1.3). 

If $(X,R)$ is an abelian monoid presentation, then we denote by $\l X\mid R\r$ the abelian monoid $\l X\r/\sim_R$ where $\sim_R$ is the congruence on $\l X\r$ generated by the relations in $R$. We say that an abelian monoid $M$ {\it is presented by $(X,R)$} if $M\cong \l X\mid R\r$.



\begin{deff}\label{defbergpres}
A {\it Bergman presentation} is a finite and good abelian monoid presentation $(X,R)$, where $R=\{a_i=b_i\mid i\in I\}$, together with a partition $I=I_{\blue}\sqcup I_{\red}$ such that conditions (i) and (ii) below are satisfied.
\begin{enumerate}[(i)]
\medskip
\item For any $i\in I_{\blue}$ there is a subset $X_i\subseteq X$ such that $|X_i|\geq 2$ and $b_i=\sum_{x\in X_i}x$ in $\l X\r$. Moreover, we require that the sets $X_i~(i\in I_{\blue})$ are pairwise disjoint. 
\bigskip
\item There is an ordering $i_1,\dots, i_m$ of the elements of $I_{\blue}$ such that $a_{i_k}\in \big\l X_{\base}\cup (\bigcup_{p=1}^{k-1} X_{i_p})\big\r$ for any $1\leq k\leq m$, where $X_{\base}=X\setminus(\bigcup_{i\in I_{\blue}} X_i)$.
\medskip
\end{enumerate}
We call $R_{\blue}=\{a_i=b_i\mid i\in I_{\blue}\}$ the {\it family of blue relations}, $R_{\red}=\{a_i=b_i\mid i\in I_{\red}\}$ the {\it family of red relations}, the elements of $X_{\base}$ {\it base generators} and an ordering of the elements of $I_{\blue}$ as in condition (ii) an {\it admissible ordering}. We call $(X,  R)$ {\it basic}, if $X=X_{\base}$ (or equivalently, if $I_{\blue}=\emptyset$).
\end{deff}

\begin{example}
Any finite and good abelian monoid presentation $(X,R)$, where $R=\{a_i=b_i\mid i\in I\}$, is a basic Bergman presentation together with the partition $I=I_{\blue}\sqcup I_{\red}$ where $I_{\blue}=\emptyset$ and $I_{\red}=I$.
\end{example}

\begin{example}\label{exbergpres}
Let $X=\{x_{0,1},x_{0,2}, x_{1,1},x_{1,2},x_{1,3},x_{2,1},x_{2,2},x_{3,1},x_{3,2}\}$, $I=\{1,2,3,4,5\}$, 
\begin{align*}
a_1&=x_{0,1}+x_{0,2},& a_2&=x_{0,1}+2x_{1,2}+x_{1,3}, & a_3&=x_{1,1}, & a_4&=x_{2,1}+x_{3,1},  &a_5&=3x_{1,3},\\
b_1&=x_{1,1}+x_{1,2}+x_{13},& b_2&=x_{2,1}+x_{2,2},& b_3&=x_{3,1}+x_{3,2},& b_4&=3x_{0,1}+x_{1,1}, &b_5&=2x_{2,2},
\end{align*}
and $R=\{a_i=b_i\mid i\in I\}$. Then $(X,R)$ is a finite and good abelian monoid presentation. Let $I_{\blue}=\{1,2,3\}$ and $I_{\red}=\{4,5\}$. Then condition (i) in Definition \ref{defbergpres} is clearly satisfied with $X_{1}=\{x_{1,1},x_{1,2},x_{1,3}\}$, $X_2=\{x_{2,1},x_{2,2}\}$ and $X_3=\{x_{3,1},x_{3,2}\}$. Condition (ii) is also satisfied, in fact there are precisely two admissible orderings of the elements of $I_{\blue}$, namely $1,2,3$ and $1,3,2$ (note that $X_{\base}=\{x_{0,1},x_{0,2}\}$). Hence $(X,R)$, together with the partition $I=I_{\blue}\sqcup I_{\red}$, is a Bergman presentation. 
\end{example}

\begin{deff}\label{defBP}
A {\it homomorphism of Bergman presentations} $\phi:(X,  R)\to (X',  R')$, where $R=\{a_i=b_i\mid i\in I\}$ and $R'=\{a'_{i'}=b'_{i'}\mid i'\in I'\}$, consists of an injective map $\phi^0:X\to X'$ and a map $\phi^1:I\to I'$ such that $\phi^1(I_{\blue})\subseteq {I'_{\blue}}$, $\psi^1(I_{\red})\subseteq I'_{\red}$, and $a'_{\phi^1(i)}=f(a_i)$ and $b'_{\phi^1(i)}=f(b_i)$ for any $i\in I$, where $f:\l X\r\to \l X'\r$ is the monoid homomorphism induced by $\phi^0$. We denote the category whose objects are the Bergman presentations and whose morphisms are the homomorphisms of Bergman presentations by ${\bf BP}$. 
\end{deff}

\subsubsection{Bergman graphs}

\begin{deff}\label{defberggraph}
A {\it Bergman graph} is a hypergraph $H$ together with a partition $H^1=H^1_{\blue}\sqcup H^1_{\red}$ such that conditions (i) and (ii) below are satisfied.
\begin{enumerate}[(i)]
\item For any $h\in H^1_{\blue}$, $r(h)$ is a set such that $|r(h)|\geq 2$. Moreover, we require that the sets $r(h)~(h\in H^1_{\blue})$ are pairwise disjoint.  
\smallskip
\item There is an ordering $h_1,\dots,h_m$ of the blue hyperedges such that $\supp(s(h_k))\subseteq H^0_{\base}\cup(\bigcup_{p=1}^{k-1}r(h_p))$ for any $1\leq k\leq m$, where $H^0_{\base}=H^0\setminus(\bigcup_{h\in H^1_{\blue}}r(h))$.
\end{enumerate}
We call the elements of $H^1_{\blue}$ {\it blue hyperedges}, the elements of $H^1_{\red}$ {\it red hyperedges}, the elements of $H^0_{\base}$ {\it base vertices} and an ordering of the blue hyperedges as in condition (ii) an {\it admissible ordering}. We call the Bergman graph $H$ {\it basic}, if $H^0=H^0_{\base}$ (or equivalently, if there are no blue hyperedges).
\end{deff}

\begin{rmk}\label{rmkberggraph}
Let $H$ be a Bergman graph that contains at least one vertex and $h_1,\dots,h_m$ an admissible ordering of the blue hyperedges. Suppose that $H^0_{\base}=\emptyset$. If $m=0$ (i.e. there are no blue hyperedges), then $H^0=H^0_{\base}=\emptyset$, a contradiction. If $m\neq 0$, then it follows from condition (ii) in Definition \ref{defberggraph} that $\supp(s(h_1))\subseteq H^0_{\base}=\emptyset$, again a contradiction. Hence $H$ contains at least one base vertex. Similarly one can show that if $( X,  R)$ is a Bergman presentation where $X\neq \emptyset$, then $X$ contains at least one base generator.
\end{rmk}

\begin{example}
Any hypergraph $H$ is a basic Bergman graph together with the partition $H^1=H^1_{\blue}\sqcup H^1_{\red}$ where $H_{\blue}^1=\emptyset$ and $H_{\red}^1=H^1$. 
\end{example}

\begin{example}\label{exberggraph}
Let $H$ be the hypergraph defined by $H^0=\{v_{0,1},v_{0,2}, v_{1,1},v_{1,2},v_{1,3},v_{2,1},v_{2,2},$ $v_{3,1},v_{3,2}\}$, $H^1=\{e,f,g,h,k\}$, $s(e)=\{v_{0,1},v_{0,2}\}$, $r(e)=\{v_{1,1},v_{1,2},v_{1,3}\}$, $s(f)=\{v_{0,1},v_{1,2},v_{1,2},v_{1,3}\}$, $r(f)=\{v_{2,1},v_{2,2}\}$, $s(g)=\{v_{1,1}\}$, $r(g)=\{v_{3,1},v_{3,2}\}$, $s(h)=\{v_{2,1},v_{3,1}\}$, $r(h)=\{v_{0,1},v_{0,1},v_{0,1},v_{1,1}\}$, $s(k)=\{v_{1,3},v_{1,3},$ $v_{1,3}\}$ and $r(k)=\{v_{2,2},v_{2,2}\}$, together with the partition $H^1=H^1_{\blue}\sqcup H^1_{\red}$ where $H_{\blue}^1=\{e,f,g\}$ and $H_{\red}^1=\{h,k\}$. We visualise $H$ as follows, where the dashed hyperedges are blue and the solid ones are red.\\
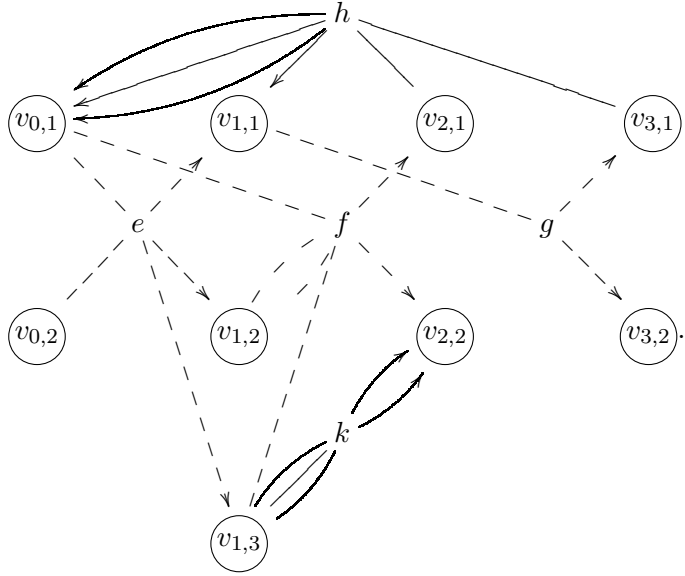
\[\xymatrix@C=18pt@R=18pt{
&&&h\ar@{-}@[red][drrr]\ar@{-}@[red][dr]&&&\\
\mathcircled{v_{0,1}}\ar@[red]@/_1pc/@{<-}[urrr]\ar@/^1pc/@[red]@{<-}[urrr]\ar@[red]@{<-}[urrr]\ar@[blue]@{--}[dr]\ar@[blue]@{--}[drrr]&&\mathcircled{v_{1,1}}\ar@[blue]@{--}[drrr]\ar@[red]@{<-}[ur]&&\mathcircled{v_{2,1}}&&\mathcircled{v_{3,1}}\\
&e\ar@[blue]@{-->}[ur]\ar@[blue]@{-->}[dr]\ar@[blue]@{-->}[dddr]&&f\ar@[blue]@{-->}[ur]\ar@[blue]@{-->}[dr]&&g\ar@[blue]@{-->}[ur]\ar@[blue]@{-->}[dr]&\\
\mathcircled{v_{0,2}}\ar@[blue]@{--}[ur]&&\mathcircled{v_{1,2}}\ar@[blue]@/^0.5pc/@{--}[ur]\ar@[blue]@/_0.5pc/@{--}[ur]&&\mathcircled{v_{2,2}}&&\mathcircled{v_{3,2}}.\\
&&&k\ar@[red]@/^0.5pc/[ur]\ar@[red]@/_0.5pc/[ur]&&&\\
&&\mathcircled{v_{1,3}}\ar@[blue]@{--}[uuur]\ar@{-}@[red]@/^0.5pc/[ur]\ar@{-}@[red]@/_0.5pc/[ur]\ar@{-}@[red][ur]&&&&
}\]
$~$\\
Condition (i) in Definition \ref{defberggraph} is clearly satisfied. Condition (ii) is also satisfied, in fact there are precisely two admissible orderings of the blue hyperedges, namely $e,f,g$ and $e,g,f$. Hence $H$ is a Bergman graph. 
\end{example}

\begin{deff}\label{defBG}
A {\it homomorphism of Bergman graphs} $\psi:H\to H'$, consists of an injective map $\psi^0:H^0\to (H')^0$ and a map $\psi^1:H^1\to (H')^1$ such that $\psi^1(H^1_{\blue})\subseteq (H')^1_{\blue}$, $\psi^1(H^1_{\red})\subseteq (H')^1_{\red}$, and 
\begin{align*}s'({\psi^1(h)})(v')&=\begin{cases}s(h)(v),\quad&\text{if }v'=\psi^0(v)\text{ for some } v\in H^0,\\
0,\quad& \text{if } v'\not\in \psi^0(H^0),\end{cases}\\
r'({\psi^1(h)})(v')&=\begin{cases}r(h)(v),\quad&\text{if }v'=\psi^0(v)\text{ for some } v\in H^0,\\
0,\quad& \text{if } v'\not\in \psi^0(H^0),\end{cases}
\end{align*}
for any $h\in H^1$ and $v'\in (H')^0$. We denote the category whose objects are the Bergman graphs and whose morphisms are the homomorphisms of Bergman graphs by ${\bf BG}$. 
\end{deff}

\subsubsection{{\bf BP} and {\bf BG} are isormorphic}
\begin{prop}\label{propbpvsbg}
The categories {\bf BP} and {\bf BG} are isomorphic. 
\end{prop}
\begin{proof}
Define a functor $\alpha:{\bf BP}\to{\bf BG}$ as follows. For an object $(X,R)$ of ${\bf BP}$, where $R=\{a_i=b_i\mid i\in I\}$, let $H=\alpha(X,R)$ be the object of {\bf BG} such that $H^0=X$, $H^1=I$, $H^1_{\blue}=I_{\blue}$, $H^1_{\red}=I_{\red}$, and for any $i\in I$ the maps $s(i):X\to \N_0$ and $r(i):X\to \N_0$ are determined by the equations $a_i=\sum_{x\in X}s(i)(x)x$ and $b_i=\sum_{x\in X}r(i)(x)x$, respectively. For a morphism $\phi:(X,R)\to (X',R')$ in ${\bf BP}$ let $\alpha(\phi):\alpha(X,R)\to \alpha(X',R')$ be the morphism such that $\alpha(\phi)^0(x)=\phi^0(x)$ for any $x\in X$ and $\alpha(\phi)^1(i)=\phi^1(i)$ of any $i\in I$.

Next define a functor $\beta:{\bf BG}\to{\bf BP}$ as follows. For an object $H$ of ${\bf BG}$ let $(X, R)=\beta(H)$, where $R=\{a_i=b_i\mid i\in I\})$, be the object of {\bf BP} such that $X=H^0$, $I=H^1$, $I_{\blue}=H^1_{\blue}$, $I_{\red}=H^1_{\red}$, and $a_h=\sum_{v\in H^0}s(h)(v)v$ and $b_h=\sum_{v\in H^0}r(h)(v)v$ for any $h\in H^1$. For a morphism $\psi:H\to H'$ in ${\bf BG}$ let $\beta(\psi):\beta(H)\to \beta(H')$ be the morphism such that $\beta(\psi)^0(v)=\psi^0(v)$ for any $v\in H^0$ and $\beta(\psi)^1(h)=\psi^1(h)$ for any $h\in H^1$.

Clearly $\alpha\circ\beta=\id_{\bf BG}$ and $\beta\circ\alpha=\id_{\bf BP}$. Thus ${\bf BP}\cong {\bf BG}$.
\end{proof}

\subsection{Bergman algebras}
\subsubsection{Bergman algebras of Bergman presentations}

\begin{deff}\label{defbergalg}
Let $(X,R)$ be a Bergman presentation where $R=\{a_i=b_i\mid i\in I\}$. Let $i_1,\dots, i_m$ be an admissible ordering of the elements of $I_{\blue}$. We will construct a $K$-algebra $B$ which has the property that $\V(B)\cong \l X\mid R\r$. We start with the algebra $A_0=K^{X_{\base}}$. Clearly $\V(A_0)\cong \l X_{\base}\r$. By applying Lemma \ref{lembergpres1.5} we obtain an algebra $A_1$ such that $\V(A_1)\cong\l X_{\base}\sqcup X_{i_1}\mid a_{i_1}=b_{i_1}\r$ (note that $a_{i_1}$ is a sum of generators from $X_{\base}$ by Definition \ref{defbergpres}(ii)). By applying Lemma \ref{lembergpres1.5} again we obtain an algebra $A_2$ such that $\V(A_2)\cong\l X_{\base}\sqcup X_{i_1}\sqcup X_{i_2}\mid a_{i_1}=b_{i_1}, a_{i_2}=b_{i_2}\r$ (note that $a_{i_1}$ is a sum of generators from $X_{\base}\sqcup X_{i_1}$ by Definition \ref{defbergpres}(ii)). We proceed like this until we obtain an algebra $A_m$ such that $\V(A_m)\cong\l X\mid a_{i}=b_{i}~(i\in I_{\blue})\r$. Now we apply Lemma \ref{lembergpres2} to obtain an algebra $B$ such that $\V(B)\cong\l X\mid R\r$ (we adjoin the remaining relations $a_{i}=b_{i}~(i\in I_{\red})$). We call the $K$-algebra $B=B(X,R)$ the {\it Bergman algebra} of $(X,R)$. 
\end{deff}

\begin{rmk}\label{remfreechoice}
Clearly, the definition of $B(X,R)$ does not depend on the chosen admissible ordering of the elements of $I_{\blue}$. In fact, one can adjoin the missing generators and relations to $\V(A_0)$ in an arbitrary order without changing the isomorphism class of $B(X,R)$, with the only restriction that one cannot adjoin a relation $a_i=b_i$ before all the generators appearing in this relation have been adjoined.
\end{rmk}

\begin{rmk}\label{rmkbergalgpres}
Let $(X,R)$ be a Bergman presentation where $R=\{a_i=b_i\mid i\in I\}$. By investigating the construction of the Bergman algebra $B=B(X,R)$ one obtains the following presentation of $B$ (in view of Lemmas \ref{lembergpres1.5} and \ref{lembergpres2}). Fix an admissible ordering $i_1,\dots,i_m$ of the elements of $I_{\blue}$. Write $X_{\base}=\{x_{0,1},\dots x_{0,t_0}\}$ and $X_{i_k}=\{x_{k,1},\dots,x_{k,t_k}\}$ for any $1\leq k\leq m$ (hence $b_{i_k}=x_{k,1}+\dots+x_{k,t_k}$ for any $1\leq k\leq m$). Moreover, write $a_{i_k}=y_{k,1}+\dots+y_{k,u_k}$ for any $1\leq k\leq m$ where $y_{k,1},\dots, y_{k,u_k}\in X_{\base}\sqcup X_{i_1}\sqcup\dots\sqcup X_{i_{k-1}}$. We will inductively associate to any element $c$ of $\l X\r\setminus\{0\}$ a quadratic matrix $\epsilon_c$ as follows. 

For any element $c\in\l X\r\setminus\{0\}$ we fix an ordering $c=x_{c,1}+\dots+x_{c,p_c}$ of its summands $x_{c,1},\dots,x_{c,p_c}\in X$. To any base generator $x_{0,l}$ where $1\leq l\leq t_0-1$ we associate a $1\times 1$ matrix $\epsilon_{x_{0,l}}$ whose only entry is a symbol. To the generator $x_{0,t_0}$ we associate the matrix $\epsilon_{x_{0,t_0}}=(1)-\epsilon_{x_{0,1}}-\dots-\epsilon_{x_{0,t_0-1}}$. Moreover, we associate to any element $c\in\l X_{\base}\r$ the matrix $\epsilon_c=\epsilon_{x_{c,1}}\oplus \dots \oplus\epsilon_{x_{c,p_c}}$.

Now suppose that for some $1\leq k\leq m-1$ we have associated to any element $c\in \l X_{\base}\sqcup X_{i_1}\sqcup\dots\sqcup X_{i_{k-1}}\r\setminus\{0\}$ a quadratic matrix $\epsilon_c$ of dimension $|\epsilon_c|$. Then we associate to any generator $x_{k,l}$ where $1\leq l \leq t_k-1$ a quadratic matrix $\epsilon_{x_{k,l}}$ of dimension $|\epsilon_{a_{i_k}}|$ whose entries are symbols. To the generator $x_{k,t_k}$ we associate the matrix $\epsilon_{x_{k,t_k}}=\epsilon_{a_{i_k}}-\epsilon_{x_{k,1}}-\dots-\epsilon_{x_{k,t_{k}-1}}$. Moreover, we associate to any element $c\in \l X_{\base}\sqcup X_{i_1}\sqcup\dots\sqcup X_{i_{k}}\r\setminus\{0\}$ the matrix $\epsilon_c=\epsilon_{x_{c,1}}\oplus \dots \oplus\epsilon_{x_{c,p_c}}$. 

Finally, we associate to any $i\in I_{\red}$ a $|\epsilon_{a_i}|\times |\epsilon_{b_i}|$ matrix $\sigma_{i}$ and a $|\epsilon_{b_i}|\times |\epsilon_{a_i}|$ matrix $\sigma'_{i}$ whose entries are symbols. With this notation the $K$-algebra $B$ has the presentation 
\begin{align*}
B\cong\Big\l& \epsilon_{x_{k,l}}~(0\leq k\leq m,~1\leq l\leq t_k-1),~\sigma_i,\sigma'_i~(i\in I_{\red})\mid\\
&\epsilon_{a_{i_k}}\epsilon_{x_{k,l}}\epsilon_{a_{i_k}}=\epsilon_{x_{k,l}}~(1\leq k\leq m,~1\leq l\leq t_k-1),\\
&\epsilon_{x_{k,l}}\epsilon_{x_{k,l'}}=\delta_{l,l'}\epsilon_{x_{k,l}}~(0\leq k\leq m,~1\leq l,l'\leq t_k-1),\\
&\epsilon_{a_i}\sigma_i\epsilon_{b_i}=\sigma_i,~\epsilon_{b_i}\sigma'_i\epsilon_{a_i}=\sigma'_i,~\sigma_i\sigma'_i=\epsilon_{a_i}, ~\sigma'_i\sigma_i=\epsilon_{b_i}~(i\in I_{\red})
\Big\r.
\end{align*}
\end{rmk}

\begin{rmk}\label{rembergfunc}
We leave it to the reader to deduce from the previous remark that $B$ defines a functor ${\bf BP}\to {\bf ALG}$ where ${\bf ALG}$ denotes the category of $K$-algebras. Hence isomorphic Bergman presentations define isomorphic Bergman algebras.
\end{rmk}

\begin{example}\label{exbergembrace}
\begin{enumerate}[(a)]
\item
It follows from Remarks \ref{remlpa}, \ref{remslpa}(b), \ref{remwlpa}(b) and \ref{remhlpa}(c) that any Leavitt algebra of a graph, separated graph, vertex-weighted graph or hypergraph is isomorphic to a Bergman algebra of a basic Bergman presentation. Conversely, any Bergman algebra of a basic Bergman presentation is isomorphic to a Leavitt path algebra of a hypergraph (follows from Remark \ref{remhlpa}(c)).
\item
It follows from Remarks \ref{remwlpa}(c) and \ref{rembslpa}(b) that any Leavitt algebra of a weighted graph or a bi-separated graph satisfying Condition (B) is isomorphic to a Bergman algebra of a not necessarily basic Bergman presentation.
\item 
Any Cohn-Leavitt path algebra based on a graph or separated graph is isomorphic to a Bergman algebra of a not necessarily basic Bergman presentation, see \cite{lpabook} and \cite[\S 4]{aragoodearl}.
\end{enumerate}
\end{example}

\begin{example}\label{exbergalgpres}
Let $(X,R)$ be the Bergman presentation from Example \ref{exbergpres}. Let $P_{0,1}$ be the finitely generated projective $K\times K$ module $K\times 0$ and $P_{0,2}$ the finitely generated projective $K\times K$ module $0\times K$. Then
\begin{align*}
B(X,R)\cong&K\times K\Big\langle [P_{1,1}],[P_{1,2}],[P_{1,3}]\Mid[P_{0,1}] +  [P_{0,2}]=[P_{1,1}] +  [P_{1,2}] +  [P_{1,3}]\Big\rangle\\
&\Big\langle [P_{2,1}], [P_{2,2}]\Mid[P_{0,1}] +  [P_{1,2}] +  [P_{1,2}] +  [P_{1,3}]=[P_{2,1}] +  [P_{2,2}]\Big\rangle\Big\langle [P_{3,1}],[P_{3,2}]\Mid[P_{1,1}]=[P_{3,1}] +  [P_{3,2}]\Big\rangle\\& \Big\langle [P_{1,3}] +  [P_{1,3}] +  [P_{1,3}]= [P_{2,2}] +  [P_{2,2}]\Big\rangle\Big\langle [P_{2,1}] +  [P_{3,1}]= [P_{1,1}] +  [P_{0,1}] +  [P_{0,1}] +  [P_{0,1}]\Big\rangle.
\end{align*}
\end{example}

\subsubsection{Bergman algebras of Bergman graphs}

\begin{deff}
Let $\beta:{\bf BG}\to {\bf BP}$ be the isomorphism of categories defined in the proof of Proposition \ref{propbpvsbg}. If $H$ is a Bergman graph, then we call the algebra $B(\beta(H))$ the {\it Bergman algebra} of $H$ and denote it by $B(H)$.
\end{deff}

\begin{rmk}\label{rembergalg}
\begin{enumerate}[(a)]
\item It follows from Remark \ref{rembergfunc} that isomorphic Bergman graphs define isomorphic Bergman algebras (since $B\circ \beta:{\bf BG}\to {\bf ALG}$ is a functor).
\medskip
\item The orientation of the red hyperedges does not play a role in the definition of $B(H)$. I.e., if $H$ and $H'$ differ only in the orientation of the red hyperedges, then $B(H)\cong B(H')$.
\medskip
\item The abelian monoid $\V(B(H))$ has the presentation 
\begin{align*}
\V(B(H))\cong\Big\l
H^0~\mid& ~\sum s(h)= \sum r(h)~(h\in H^1)\Big\r
\end{align*}
where $\sum s(h)=\sum_{v\in H^0}s(h)(v)v$ and $\sum r(h)=\sum_{v\in H^0}r(h)(v)v$.
\medskip
\item If $H$ is basic, then $B(H)\cong L(H)$ where $L(H)$ is the Leavitt path algebra of $H$ (follows from Remark \ref{remhlpa}(c)).
\end{enumerate}
\end{rmk}

\begin{example}
By Example \ref{exbergembrace}, any Leavitt path algebra of a graph, separated graph, vertex-weighted graph or hypergraph is isomorphic to a Bergman algebra of a basic Bergman graph. Moreover, any Leavitt algebra of a weighted graph or a bi-separated graph satisfying Condition (B), as well as any Cohn-Leavitt path algebra based on a graph or separated graph is isomorphic to a Bergman algebra of a not necessarily basic Bergman graph.
\end{example}

\begin{example}\label{exbergalglpa}
Leavitt path algebras of graphs can be realised as Bergman algebras of Bergman graphs as follows. Let $E=(E^0,E^1, s,r)$ be a graph. Let $H=(H^0, H^1, s_H, r_H)$ be the basic Bergman graph defined by $H^0=E^0$, $H^1=E^{0}_{\reg}$, $s_H(v)(w)=\delta_{v,w}v$ and $r_H(v)(w)=\#\{e\in s^{-1}(v)\mid r(e)=w\}$ for any $v\in H^1$ and $w\in H^0$. Then $L(E)\cong B(H)$ (follows from Remark \ref{remlpa}).
\end{example}

\begin{example}\label{exbergalggraph2}
Let $H$ be the basic Bergman graph\\
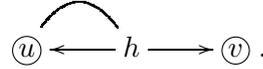
\[\xymatrix{
\mathcircled{u}\ar@{-}@[red]@/^1.5pc/[r]&h\ar@{->}@[red][l]\ar@{->}@[red][r]&\mathcircled{v}
}.\]
Then $B(H)=K\times K\big\langle [K\times 0]= [K\times K]\big\rangle$. It follows from Remark \ref{remlpa} (or Example \ref{exbergalglpa}) that $B(H)$ is isomorphic to the Leavitt path algebra of the Toeplitz graph 
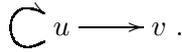
\[\xymatrix{u\ar@(dl,ul)\ar[r]&v}.\]
\end{example}

\begin{example}\label{exbergalggraph3}
Let $H$ be the basic Bergman graph\\
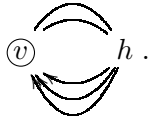
\[\xymatrix{
\mathcircled{v}\ar@{-}@[red]@/^1.5pc/[r]\ar@{-}@[red]@/^1pc/[r]&h\ar@/^1pc/@{->}@[red][l]\ar@/^1.5pc/@{->}@[red][l]\ar@/^2pc/@{->}@[red][l]
}.\]
$~$\\
\\
Then $B(H)=K\big\langle [K^2]=[K^3]\big\rangle$. It follows from Remark \ref{remwlpa}(b) that $B(H)$ is isomorphic to the Leavitt path algebra of the vertex-weighted graph 
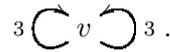
\[\hspace{0.4cm}\xymatrix{v\ar@(dl,ul)^{3}\ar@(dr,ur)_{3}}.\]
\end{example}

\begin{example}\label{exbergalggraph5}
Let $H$ be the basic Bergman graph\\
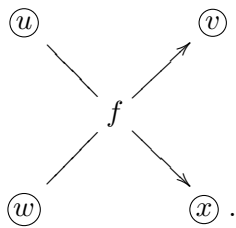
\[\xymatrix@C=18pt@R=18pt{
\mathcircled{u}\ar@[red]@{-}[dr]&&\mathcircled{v}\\
&f\ar@[red]@{->}[ur]\ar@[red]@{->}[dr]&\\
\mathcircled{w}\ar@[red]@{-}[ur]&&\mathcircled{x}~.
}\]
$~$\\
\\
Then $B(H)=K\times K\times K\times K\big\langle [K\times K\times 0\times 0]=[0\times 0\times K\times K]\big\rangle$. By Remark \ref{rembergalg}(d) $B(H)\cong L(H)$.
\end{example}

\begin{example}\label{exbergalggraph4}
Let $H$ be the Bergman graph\\
\[\xymatrix@C=18pt@R=18pt{
&h\ar@{->}@[red][drrr]&&&\\
&h'\ar@{->}@[red][r]&\mathcircled{v_{1,1}}&&\mathcircled{v_{2,1}}\\
\mathcircled{v_{0,1}}\ar@[red]@{-}[ur]\ar@[red]@{-}[dr]\ar@/^1.3pc/@[red]@{-}[dr]\ar@/_1.2pc/@[red]@{-}[dr]\ar@/^1.2pc/@[red]@{-}[uur]\ar@[red]@{-}[uur]\ar@[blue]@{--}[r]\ar@/^1.2pc/@[blue]@{--}[rrr]&e\ar@[blue]@{-->}[ur]\ar@[blue]@{-->}[dr]&&f\ar@[blue]@{-->}[ur]\ar@[blue]@{-->}[dr]&\\
&h''\ar@/_2pc/@{->}@[red][rrr]\ar@/^3pc/@{->}@[red][ul]&\mathcircled{v_{1,2}}\ar@[blue]@{--}[ur]&&\mathcircled{v_{2,2}}~.
}\]
$~$\\
\\
Then 
\begin{align*}
B(H)=&K\big\langle [P_{1,1}],[P_{1,2}]\Mid [P_{0,1}]=[P_{1,1}] +  [P_{1,2}]\big\rangle\big\langle [P_{2,1}],[P_{2,2}]\Mid [P_{0,1}] +  [P_{1,2}]=[P_{2,1}] +  [P_{2,2}]\big\rangle\\&\big\langle [P_{0,1}] +  [P_{0,1}]= [P_{2,1}]\big\rangle\big\langle [P_{0,1}]= [P_{1,1}]\big\rangle\big\langle [P_{0,1}] +  [P_{0,1}] +  [P_{0,1}]= [P_{0,1}] +  [P_{2,2}]\big\rangle
\end{align*}
where $P_{0,1}$ is the free $K$-module of rank $1$. It follows from \cite[\S 4]{preusser-1} that $B(H)$ is isomorphic to the Leavitt path algebra of the weighted graph 
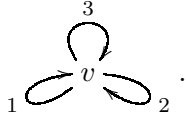
\[\hspace{0.4cm}\xymatrix{v\ar@(dl,l)^{1}\ar@(r,dr)^{2}\ar@(ul,ur)^{3}}.\]
\end{example}

\begin{example}\label{exbergalggraph}
Let $H$ be the Bergman graph from Example \ref{exberggraph}. Clearly $\beta(H)\cong(X,R)$ where $(X,R)$ is the Bergman presentation from Example \ref{exbergpres}. Hence, in view of Example \ref{exbergalgpres},
\begin{align*}
B(H)\cong&K\times K\Big\langle [P_{1,1}],[P_{1,2}],[P_{1,3}]\Mid[P_{0,1}] +  [P_{0,2}]=[P_{1,1}] +  [P_{1,2}] +  [P_{1,3}]\Big\rangle\\
&\Big\langle [P_{2,1}], [P_{2,2}]\Mid[P_{0,1}] +  [P_{1,2}] +  [P_{1,2}] +  [P_{1,3}]=[P_{2,1}] +  [P_{2,2}]\Big\rangle\Big\langle [P_{3,1}],[P_{3,2}]\Mid[P_{1,1}]=[P_{3,1}] +  [P_{3,2}]\Big\rangle\\& \Big\langle [P_{1,3}] +  [P_{1,3}] +  [P_{1,3}]= [P_{2,2}] +  [P_{2,2}]\Big\rangle\Big\langle [P_{2,1}] +  [P_{3,1}]= [P_{1,1}] +  [P_{0,1}] +  [P_{0,1}] +  [P_{0,1}]\Big\rangle.
\end{align*}
where $P_{0,1}=K\times 0$ and $P_{0,2}=0\times K$.
\end{example}

\section{Moves for Bergman algebras}
In this section we define the following moves for Bergman presentations and Bergman graphs and prove that they preserve the isomorphism class of the corresponding Bergman algebras: red shift move, blue shift move, enqueuing and outsplitting. Moreover, we define the following moves for basic Bergman presentations and basic Bergman graphs and prove that they preserve the Morita equivalence class of the corresponding Bergman algebras: lonely generator/vertex elimination, collapsing and insplitting. We leave it to the reader to recover (via Example \ref{exbergalglpa}) the results, that in the context of usual graphs the shift move and outsplitting preserve the isomorphism class of the corresponding Leavitt path algebra, while source elimination, collapsing and insplitting preserve the Morita equivalence class of the corresponding Leavitt path algebra, see \cite[Theorem 1.1]{kbgm}, \cite[\S 6.2, \S 6.3]{lpabook} and \cite[\S 2]{AN}. 

We also show that any collapsing is a composition of a finite number of red shift moves and lonely generator/vertex eliminations, and that any insplitting is a composition of a finite number of red shift moves and extensions, where by an extension we mean the move that is inverse to a lonely generator/vertex elimination. These two results seem to be new even in the context of usual graphs (a red shift move corresponds to a shift move and a lonely vertex elimination corresponds to a source elimination in this context).  

Throughout this section $(X,R)$ denotes a Bergman presentation, where $R=\{a_i=b_i\mid i\in I\}$. We fix an admissible ordering $i_1,\dots, i_m$ of the elements of $I_{\blue}$. 
Moverover, we denote the Bergman graph corresponding to $(X,R)$ by $H$ (see \S 3.1.3). Clearly the admissible ordering $i_1,\dots, i_m$ of the elements of $I$ induces an admissible ordering $h_1,\dots, h_m$ of the blue hyperedges in $H$. 

\subsection{Red shift move}
\begin{deff}
Let $i\in I_{\red}$ and $a,b\in  \l X\r$ such that $a_{i}=a$ and $b_i=b$ in the abelian monoid $\l X\mid a_{j}=b_{j}~(j\in I\setminus\{i\})\r$. Let $(X,R')$ be the Bergman presentation obtained from $(X,R)$ by replacing the $i$-th relation $a_{i}=b_{i}$ by the relation $a=b$ (of the same colour, namely red). Then we say that $(X,R')$ can be obtained from $(X,R)$ by a {\it red shift move}.
\end{deff}

\begin{thm}\label{thmshiftredp}
Suppose the Bergman presentation $(X,R')$ can be obtained from $(X,R)$ by a red shift move. Then $B(X,R')\cong B(X,R)$.
\end{thm}
\begin{proof}
The theorem follows Definition \ref{defbergalg} and Remark \ref{remfreechoice} (in the construction of $B(X,R)$ (respectively $B(X,R')$) adjoin the relation $a_{i}=b_{i}$ (resp. $a=b$) after having adjoined all generators and all other relations).
\end{proof}

\begin{deff}
Let $h\in H^1_{\red}$ and $u_1,\dots,u_p,v_1,\dots,v_q\in H^0$ such that $\sum s(h)=\sum_{i=1}^pu_i$ and $\sum r(h)=\sum_{i=1}^qv_i$ in the abelian monoid $\big\l H^0\mid \sum s(g)=\sum r(g)~(g\in H^1\setminus \{h\})\big\r$. Let $H'$ be the the Bergman graph obtained from $H$ by replacing the hyperedge $h$ by a red hyperedge $h'$ whose source and range multisets are defined by $s(h')=\{u_1,\dots,u_p\}$ and $r(h')=\{v_1,\dots,v_q\}$, respectively. Then we say that $H'$ can be obtained from $H$ by a {\it red shift move}.
\end{deff}

\begin{thm}\label{thmshiftred}
Suppose the Bergman graph $H'$ can be obtained from $H$ by a red shift move. Then $B(H')\cong B(H)$. 
\end{thm}
\begin{proof}
The theorem follows from Theorem \ref{thmshiftredp}.
\end{proof}

\begin{example}
Let $H$ and $H'$ be the (basic) Bergman graphs
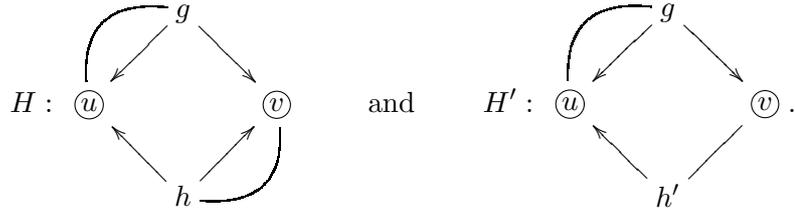
\[H:\,\vcenter{\vbox{\xymatrix@R=20pt@C=20pt{
&g\ar@[red]@{->}[dr]\ar@[red]@{->}[dl]&\\
\mathcircled{u}\ar@/^1.5pc/@[red]@{-}[ur]&&\mathcircled{v}\ar@/^1.5pc/@[red]@{-}[dl]\\
&h\ar@[red]@{->}[ur]\ar@[red]@{->}[ul]&
}}}
\quad\quad \text{ and }\quad\quad H':\,\vcenter{\vbox{\xymatrix@R=20pt@C=20pt{
&g\ar@[red]@{->}[dr]\ar@[red]@{->}[dl]&\\
\mathcircled{u}\ar@/^1.5pc/@[red]@{-}[ur]&&\mathcircled{v}\\
&h'\ar@[red]@{-}[ur]\ar@[red]@{->}[ul]&
}}}.\]
Then $B(H)\cong B(H')$ by Theorem \ref{thmshiftred} since $u+v=u$ in the abelian monoid $\langle u,v\mid u=u+v\rangle$.
\end{example}

\begin{example}
Let $H$ and $H'$ be the Bergman graphs\\
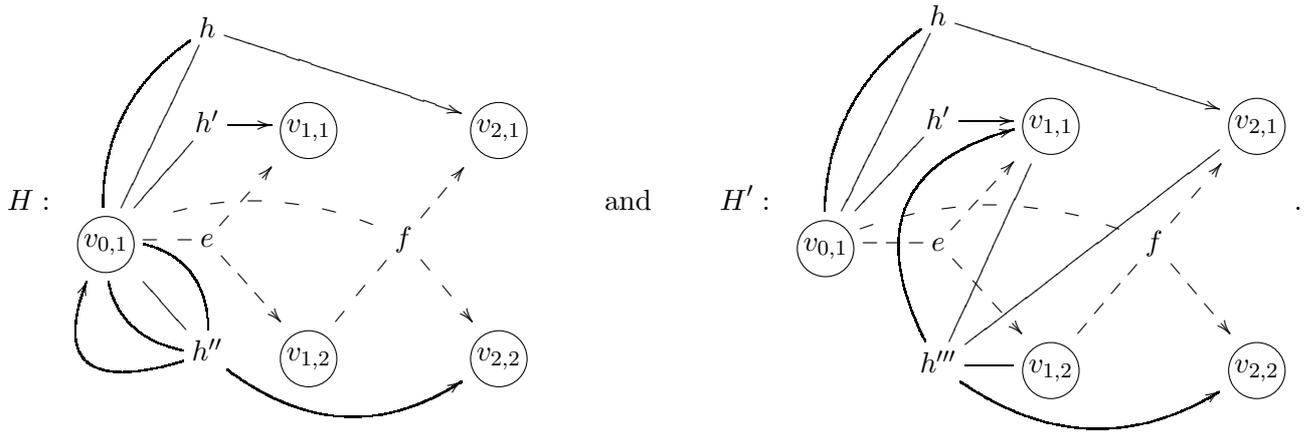
\[H:\,\,\vcenter{\vbox{\xymatrix@C=15pt@R=15pt{
&h\ar@{->}@[red][drrr]&&&\\
&h'\ar@{->}@[red][r]&\mathcircled{v_{1,1}}&&\mathcircled{v_{2,1}}\\
\mathcircled{v_{0,1}}\ar@[red]@{-}[ur]\ar@[red]@{-}[dr]\ar@/^1.3pc/@[red]@{-}[dr]\ar@/_1.2pc/@[red]@{-}[dr]\ar@/^1.2pc/@[red]@{-}[uur]\ar@[red]@{-}[uur]\ar@[blue]@{--}[r]\ar@/^1.2pc/@[blue]@{--}[rrr]&e\ar@[blue]@{-->}[ur]\ar@[blue]@{-->}[dr]&&f\ar@[blue]@{-->}[ur]\ar@[blue]@{-->}[dr]&\\
&h''\ar@/_2pc/@{->}@[red][rrr]\ar@/^3pc/@{->}@[red][ul]&\mathcircled{v_{1,2}}\ar@[blue]@{--}[ur]&&\mathcircled{v_{2,2}}
}}}
\quad\quad\text{ and }\quad \quad H':\,\,\vcenter{\vbox{\xymatrix@C=18pt@R=18pt{
&h\ar@{->}@[red][drrr]&&&\\
&h'\ar@{->}@[red][r]&\mathcircled{v_{1,1}}&&\mathcircled{v_{2,1}}\\
\mathcircled{v_{0,1}}\ar@[red]@{-}[ur]\ar@/^1.2pc/@[red]@{-}[uur]\ar@[red]@{-}[uur]\ar@[blue]@{--}[r]\ar@/^1.2pc/@[blue]@{--}[rrr]&e\ar@[blue]@{-->}[ur]\ar@[blue]@{-->}[dr]&&f\ar@[blue]@{--}[dl]\ar@[blue]@{-->}[ur]\ar@[blue]@{-->}[dr]&\\
&h'''\ar@{-}@[red][r]\ar@{-}@[red][uur]\ar@{-}@[red][uurrr]\ar@/_2pc/@{->}@[red][rrr]\ar@/^3pc/@{->}@[red][uur]&\mathcircled{v_{1,2}}&&\mathcircled{v_{2,2}}
}}}.\]
$~$\\
\\
Then $B(H)\cong B(H')$ by Theorem \ref{thmshiftred} since $v_{0,1}+v_{0,1}+v_{0,1}=v_{1,1}+v_{1,2}+v_{2,1}$ and $v_{0,1}+v_{2,2}=v_{1,1}+v_{2,2}$ in the abelian monoid $\big\l H^0\mid \sum s(g)=\sum r(g)~(g\in H^1\setminus \{h''\})\big\r$.
\end{example}

\subsection{Blue shift move }

\begin{deff}
Let $k\in\{1,\dots,m\}$ and $a\in  \l X_{\base}\sqcup X_{i_1}\sqcup\dots\sqcup X_{i_{k-1}}\r$ such that $a_{i_k}=a$ in the abelian monoid $\big\l X_{\base}\sqcup X_{i_1}\sqcup\dots\sqcup X_{i_{k-1}}\mid a_{i_p}=b_{i_p}~(1\leq p\leq k-1),~a_{j}=b_{j}~(j\in J)\big\r$, where $J\subseteq I_{\red}$ consists of all indices $j$ such that $a_j,b_j\in  \l X_{\base}\sqcup X_{i_1}\sqcup\dots\sqcup X_{i_{k-1}}\r$. Let $(X,R')$ be the Bergman presentation obtained from $(X,R)$ by replacing the $i_k$-th relation $a_{i_k}=b_{i_k}$ by the relation $a=b_{i_k}$ (of the same colour, namely blue). Then we say that $(X,R')$ can be obtained from $(X,R)$ by a {\it blue shift move}.
\end{deff}

\begin{thm}\label{thmshiftbluep}
Suppose the Bergman presetation $(X,R')$ can be obtained from $(X,R)$ by a blue shift move. Then $B(X,R')\cong B(X,R)$.
\end{thm}
\begin{proof}
The theorem follows Definition \ref{defbergalg} and Remark \ref{remfreechoice} (in the construction of $B(X,R)$ (respectively $B(X,R')$) adjoin first the generators from $X_{\base}\sqcup X_{i_1}\sqcup\dots\sqcup X_{i_{k-1}}$ and the relations $a_{i_p}=b_{i_p}~(1\leq p\leq k-1)$ to $\V(K^{X_{\base}})$, then the relations $a_{j}=b_{j}~(j\in J)$).
\end{proof}

\begin{deff}
Let $k\in\{1,\dots,m\}$ and $u_1,\dots,u_q\in H^0_{\base}\sqcup(\bigsqcup_{p=1}^{k-1}r(h_p))$ such that $\sum s(h_k)=\sum_{i=1}^qu_i$ in the abelian monoid 
\[\Big\l H^0_{\base}\sqcup(\bigsqcup_{p=1}^{k-1}r(h_p))\Mid \sum s(h_{p})=\sum r(h_{p})~(1\leq p\leq k-1),~\sum s(g)=\sum r(g)~(g\in G)\Big\r\]
where $G\subseteq H^1_{\red}$ consists of all red hyperedges $g$ such that all sources and ranges of $g$ are contained in $H^0_{\base}\sqcup(\bigsqcup_{p=1}^{k-1}r(h_p))$. Let $H'$ be the the Bergman graph obtained from $H$ by replacing the hyperedge $h_k$ by a blue hyperedge $h$ whose source and range multisets are defined by $s(h)=\{u_1,\dots,u_q\}$ and $r(h)=r(h_k)$, respectively. Then we say that $H'$ can be obtained from $H$ by a {\it blue shift move}.
\end{deff}

\begin{thm}\label{thmshiftblue}
Suppose the Bergman graph $H'$ can be obtained from $H$ by a blue shift move. Then $B(H')\cong B(H)$. 
\end{thm}
\begin{proof}
The theorem follows from Theorem \ref{thmshiftbluep}.
\end{proof}

\begin{example}
Let $H$ and $H'$ be the Bergman graphs\\
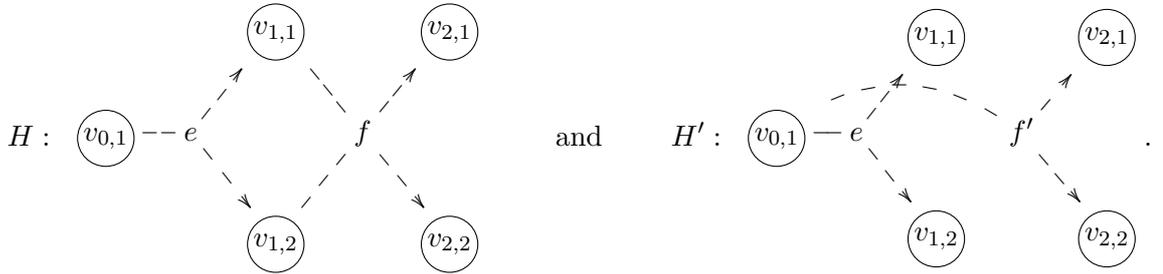
\[H:\,\,\vcenter{\vbox{\xymatrix@C=12pt@R=12pt{
&&\mathcircled{v_{1,1}}\ar@[blue]@{--}[dr]&&\mathcircled{v_{2,1}}\\
\mathcircled{v_{0,1}}\ar@[blue]@{--}[r]&e\ar@[blue]@{-->}[ur]\ar@[blue]@{-->}[dr]&&f\ar@[blue]@{-->}[ur]\ar@[blue]@{-->}[dr]&\\
&&\mathcircled{v_{1,2}}\ar@[blue]@{--}[ur]&&\mathcircled{v_{2,2}}
}}}\quad\quad\text{ and }\quad \quad H':\,\,\vcenter{\vbox{\xymatrix@C=10pt@R=10pt{
&&\mathcircled{v_{1,1}}&&\mathcircled{v_{2,1}}\\
\mathcircled{v_{0,1}}\ar@[blue]@/^1.5pc/@{--}[rrr]\ar@[blue]@{--}[r]&e\ar@[blue]@{-->}[ur]\ar@[blue]@{-->}[dr]&&f'\ar@[blue]@{-->}[ur]\ar@[blue]@{-->}[dr]&\\
&&\mathcircled{v_{1,2}}&&\mathcircled{v_{2,2}}
}}}.\]
Then $B(H)\cong B(H')$ by Theorem \ref{thmshiftblue} since $v_{0,1}=v_{1,1}+v_{1,2}$ in the abelian monoid $\langle v_{0,1},v_{1,1},v_{1,2}\mid v_{0,1}=v_{1,1}+v_{1,2}\rangle$.
\end{example}

\begin{example}
Let $H$ and $H'$ be the Bergman graphs\\
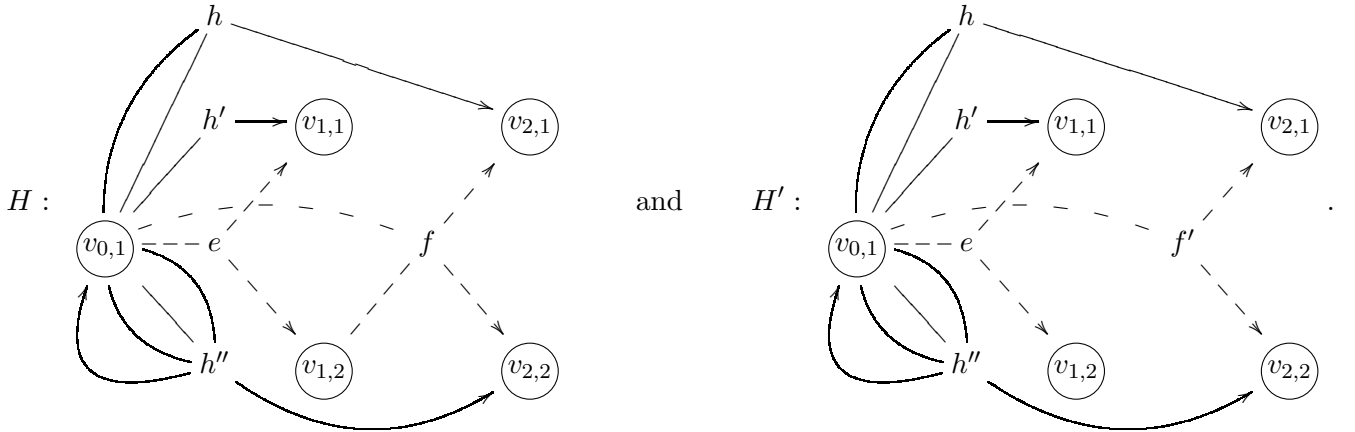
\[H:\,\,\vcenter{\vbox{\xymatrix@C=18pt@R=18pt{
&h\ar@{->}@[red][drrr]&&&\\
&h'\ar@{->}@[red][r]&\mathcircled{v_{1,1}}&&\mathcircled{v_{2,1}}\\
\mathcircled{v_{0,1}}\ar@[red]@{-}[ur]\ar@[red]@{-}[dr]\ar@/^1.3pc/@[red]@{-}[dr]\ar@/_1.2pc/@[red]@{-}[dr]\ar@/^1.2pc/@[red]@{-}[uur]\ar@[red]@{-}[uur]\ar@[blue]@{--}[r]\ar@/^1.2pc/@[blue]@{--}[rrr]&e\ar@[blue]@{-->}[ur]\ar@[blue]@{-->}[dr]&&f\ar@[blue]@{-->}[ur]\ar@[blue]@{-->}[dr]&\\
&h''\ar@/_2pc/@{->}@[red][rrr]\ar@/^3pc/@{->}@[red][ul]&\mathcircled{v_{1,2}}\ar@[blue]@{--}[ur]&&\mathcircled{v_{2,2}}
}}}
\quad\quad\text{ and }\quad\quad
H':\,\,\vcenter{\vbox{\xymatrix@C=18pt@R=18pt{
&h\ar@{->}@[red][drrr]&&&\\
&h'\ar@{->}@[red][r]&\mathcircled{v_{1,1}}&&\mathcircled{v_{2,1}}\\
\mathcircled{v_{0,1}}\ar@[red]@{-}[ur]\ar@[red]@{-}[dr]\ar@/^1.3pc/@[red]@{-}[dr]\ar@/_1.2pc/@[red]@{-}[dr]\ar@/^1.2pc/@[red]@{-}[uur]\ar@[red]@{-}[uur]\ar@[blue]@{--}[r]\ar@/^1.2pc/@[blue]@{--}[rrr]&e\ar@[blue]@{-->}[ur]\ar@[blue]@{-->}[dr]&&f'\ar@[blue]@{-->}[ur]\ar@[blue]@{-->}[dr]&\\
&h''\ar@/_2pc/@{->}@[red][rrr]\ar@/^3pc/@{->}@[red][ul]&\mathcircled{v_{1,2}}&&\mathcircled{v_{2,2}}
}}}.\]
$~$\\
\\
Then $B(H)\cong B(H')$ by Theorem \ref{thmshiftblue} since $v_{0,1}+v_{1,2}=v_{0,1}$ in the abelian monoid $\langle v_{0,1},v_{1,1},v_{1,2}\mid v_{0,1}=v_{1,1}+v_{1,2},~v_{0,1}=v_{1,1}\rangle$.
\end{example}

\subsection{Enqueuing}
\begin{deff}
Let $i\in I_{\blue}$ such that $a_{i}\in X$. For any $j\in I\setminus\{i\}$ let $a'_j$ (resp. $b'_j$) be the element of $\l X\r $ obtained from $a_{j}$ (resp. $b_{j}$) by replacing in the expression of $a_j$ (resp. $b_j$) as a sum of generators any occurrence of $a_{i}$ by $b_{i}$. Let $(X',R')$ be the Bergman presentation obtained from $(X,R)$ by removing the $i$-th relation $a_{i}=b_{i}$ and the generator $a_{i}$, and replacing each of the other relations $a_j=b_j$, where $j\in I\setminus\{i\}$, by the corresponding relation $a'_j=b'_j$ (of the same colour). Then we say that $(X',R')$ can be obtained from $(X,R)$ by {\it enqueuing}.
\end{deff}
\begin{thm}\label{thmenqp}
Suppose the Bergman presentation $(X',R')$ can be obtained from $(X,R)$ by enqueuing. Then $B( X',R')\cong B(X,R)$.
\end{thm}
\begin{proof}
The theorem follows Definition \ref{defbergalg}, Remark \ref{remfreechoice}, Theorem \ref{thmshiftredp}, Theorem \ref{thmshiftbluep} and the definition of the algebras $R\Big\l [P_1],\dots [P_t]\Mid [P]=[P_1] + \dots  +  [P_t]\Big\r$.
\end{proof}

\begin{deff}
Let $h\in H_{\blue}^1$ such that $s(h)=\{v\}$ for some $v\in H^0$. Let $H'$ be the the Bergman graph obtained from $H$ by removing the hyperedge $h$ and the vertex $v$, and replacing $s$ and $r$ by $s'$ and $r'$, respectively, which are defined by
\begin{align*}
s'(g)(u)=\begin{cases}
s(g)(u),\quad&\text{if } u\not\in r(h),\\
s(g)(u)+s(g)(v),\quad&\text{if }u\in r(h),
\end{cases}\quad\text{ and }\quad r'(h)(u)=\begin{cases}
r(g)(u),\quad&\text{if } u\not\in r(h),\\
r(g)(u)+r(g)(v),\quad&\text{if }u\in r(h).
\end{cases}
\end{align*}
Then we say that $H'$ can be obtained from $H$ by {\it enqueuing}.
\end{deff}

\begin{thm}\label{thmenq}
Suppose the Bergman graph $H'$ can be obtained from $H$ by enqueuing. Then $B(H')\cong B(H)$. 
\end{thm}
\begin{proof}
The theorem follows from Theorem \ref{thmenqp}.
\end{proof}

\begin{example}
Let $H$ and $H'$ be the Bergman graphs\\
\[H:\,\,\vcenter{\vbox{\xymatrix@C=15pt@R=15pt{
\mathcircled{v_{1,1}}\ar@[blue]@{--}[dr]&&\mathcircled{v_{2,1}}\\
&h\ar@[blue]@{-->}[ur]\ar@[blue]@{-->}[dr]&\\
\mathcircled{v_{1,2}}&&\mathcircled{v_{2,2}}
}}}\quad\quad\text{ and }\quad \quad H':\,\,\vcenter{\vbox{\xymatrix@C=15pt@R=15pt{
\mathcircled{v_{2,1}}\\
\mathcircled{v_{2,1}}\\
\mathcircled{v_{1,2}}
}}}.\]
Then $B(H)\cong B(H')$ by Theorem \ref{thmenq}.
\end{example}

\begin{example}
Let $H$ and $H'$ be the Bergman graphs\\
\[H:\,\,\vcenter{\vbox{\xymatrix@C=18pt@R=18pt{
&h'\ar@{->}@[red][r]&\mathcircled{v_{1,1}}&&\mathcircled{v_{2,1}}\\
\mathcircled{v_{0,1}}\ar@[red]@{-}[ur]\ar@[blue]@{--}[r]\ar@/^1.2pc/@[blue]@{--}[rrr]&e\ar@[blue]@{-->}[ur]\ar@[blue]@{-->}[dr]&&f\ar@[blue]@{-->}[ur]\ar@[blue]@{-->}[dr]&\\
&&\mathcircled{v_{1,2}}\ar@[blue]@{--}[ur]&&\mathcircled{v_{2,2}}
}}}\quad\quad\text{ and }\quad\quad H':\,\,\vcenter{\vbox{\xymatrix@C=18pt@R=18pt{
&\mathcircled{v_{1,1}}\ar@[blue]@{--}[dr]&&\mathcircled{v_{2,1}}\\
h'\ar@/^1.5pc/@{->}@[red][ur]\ar@{-}@[red][ur]\ar@{-}@[red][dr]&&f\ar@[blue]@{-->}[ur]\ar@[blue]@{-->}[dr]&\\
&\mathcircled{v_{1,2}}\ar@[blue]@/^1.2pc/@{--}[ur]\ar@[blue]@/_1.2pc/@{--}[ur]&&\mathcircled{v_{2,2}}
}}}.\]
Then $B(H)\cong B(H')$ by Theorem \ref{thmenq}.
\end{example}

\subsection{Outsplitting}

\begin{deff}
Let $i\in I_{\red}$ and $c_1,\dots,c_t\in \l X\r\setminus\{0\}$, where $t\geq 2$, such that $a_{i}=c_1+\dots+c_t$ in $\l X\r$. Let $(X',R')$ be the Bergman presentation obtained from $(X,R)$ by removing the $i$-th relation $a_{i}=b_{i}$ and adding generators $x_1,\dots,x_t$ to $X$, one blue relation $b_i=x_1+\dots+x_t$ and red relations $c_1=x_1,\dots, c_t=x_t$. Then we say that $(X',R')$ can be obtained from $(X,R)$ by {\it outsplitting}.
\end{deff}

\begin{thm}\label{thmoutsplitp}
Suppose the Bergman presentation $(X',R')$ can be obtained from $(X,R)$ by outsplitting. Then $B(X',R')\cong B(X,R)$.
\end{thm}
\begin{proof}
The theorem follows from Definition \ref{defbergalg} and Lemma \ref{lemoutsplit}.
\end{proof}

\begin{deff}\label{defoutsplit}
Let $h\in H^1_{\red}$ and $m_1,\dots, m_t$, where $t\geq 2$, be nonempty multisets such that $m_1+\dots+m_t=s(h)$. Let $H'$ be the the Bergman graph obtained from $H$ by removing the hyperedge $h$ and adding $t$ vertices $v_1,\dots, v_t$, one blue hyperedge $g$ such that $s(g)=r(h)$ and $r(g)=\{v_1,\dots,v_t\}$, and red hyperedges $h_1,\dots, h_t$ such that $s(h_i)=m_i$ and $r(h_i)=\{v_i\}$ for any $1\leq i\leq t$. Then we say that $H'$ can be obtained from $H$ by {\it outsplitting}.
\end{deff}

\begin{thm}\label{thmoutsplit}
Suppose the Bergman graph $H'$ can be obtained from $H$ by outsplitting. Then $B(H')\cong B(H)$. 
\end{thm}
\begin{proof}
The theorem follows from Theorem \ref{thmoutsplitp}.
\end{proof}

\begin{example}
Let $H$ and $H'$ be the Bergman graphs\\
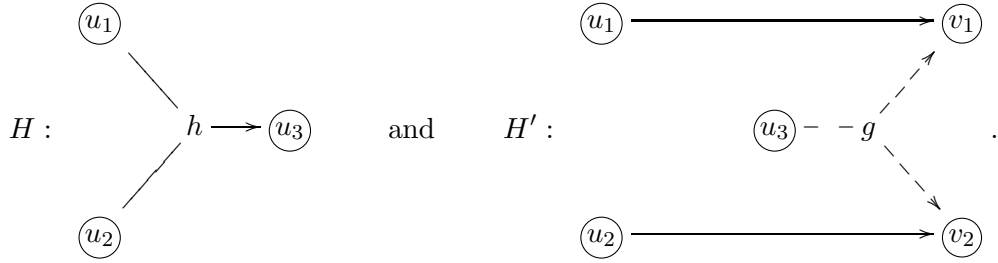
\[H:\,\,\vcenter{\vbox{\xymatrix@C=18pt@R=18pt{
\mathcircled{u_1}&&\\
&h\ar@[red]@{-}[ul]\ar@[red]@{-}[dl]&\mathcircled{u_3}\ar@[red]@{<-}[l]\\
\mathcircled{u_2}&&
}}}\quad\quad\text{ and }\quad\quad H':\,\,\vcenter{\vbox{\xymatrix@C=18pt@R=18pt{
\mathcircled{u_1}\ar@[red]@{->}[rrrr]&&&&\mathcircled{v_1}\\
&&\mathcircled{u_3}\ar@[blue]@{--}[r]&g\ar@[blue]@{-->}[ur]\ar@[blue]@{-->}[dr]&\\
\mathcircled{u_2}\ar@[red]@{->}[rrrr]&&&&\mathcircled{v_2}
}}}.\]
Then $B(H)\cong B(H')$ by Theorem \ref{thmoutsplit}. Let $H''$ be the Bergman graph 
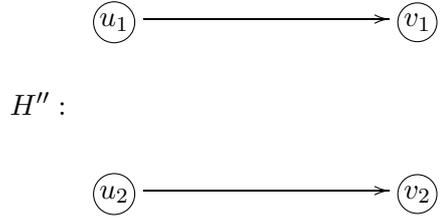
\[H'':\,\,\vcenter{\vbox{\xymatrix@C=18pt@R=18pt{
\mathcircled{u_1}\ar@[red]@{->}[rrrr]&&&&\mathcircled{v_1}\\
&&&\\
\mathcircled{u_2}\ar@[red]@{->}[rrrr]&&&&\mathcircled{v_2}
}}}.\]
Then $B(H')\cong B(H'')$ by Theorem \ref{thmenq}. In fact, if the hyperedge $h$ in Definition \ref{defoutsplit} has the property that $r(h)=\{u\}$ for some vertex $u$, then one can always apply enqueuing after having applied outsplitting.
\end{example}

\begin{example}
Let $H$ and $H'$ be the Bergman graphs\\
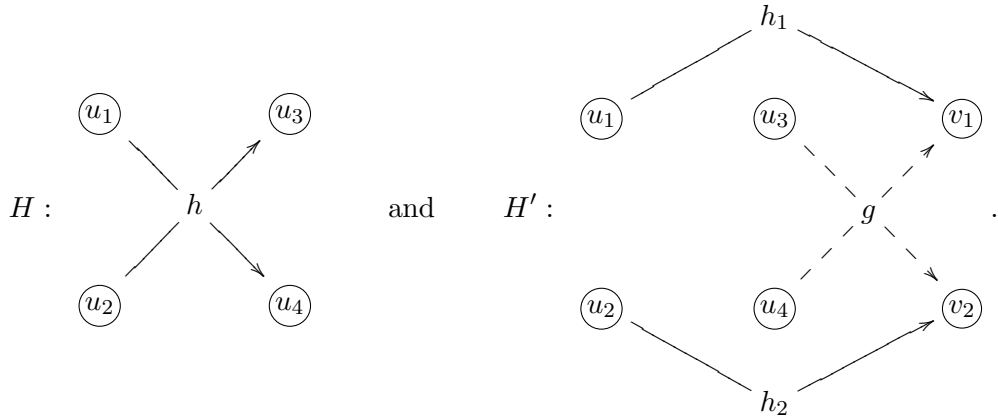
\[H:\,\,\vcenter{\vbox{\xymatrix@C=18pt@R=18pt{
\mathcircled{u_1}&&\mathcircled{u_3}\ar@[red]@{<-}[dl]\\
&h\ar@[red]@{-}[ul]\ar@[red]@{-}[dl]&\\
\mathcircled{u_2}&&\mathcircled{u_4}\ar@[red]@{<-}[ul]
}}}\quad\quad\text{ and }\quad\quad H':\,\,\vcenter{\vbox{\xymatrix@C=18pt@R=18pt{
&&h_1\ar@[red]@{->}[drr]&&\\
\mathcircled{u_1}\ar@[red]@{-}[urr]&&\mathcircled{u_3}\ar@[blue]@{--}[dr]&&\mathcircled{v_1}\\
&&&g\ar@[blue]@{-->}[ur]\ar@[blue]@{-->}[dr]&\\
\mathcircled{u_2}\ar@[red]@{-}[drr]&&\mathcircled{u_4}\ar@[blue]@{--}[ur]&&\mathcircled{v_2}\\
&&h_2\ar@[red]@{->}[urr]&&
}}}.\]
Then $B(H)\cong B(H')$ by Theorem \ref{thmoutsplit}. This time we cannot apply enqueuing to $H'$ since $s(g)=r(h)$ does not equal $\{u\}$ for some vertex $u$.
\end{example}

\subsection{Lonely vertex elimination and lonely generator elimination}
\begin{deff}\label{deflone}
We call a vertex $v\in H^0$ {\it lonely} if there is a hyperedge $h$ such that $s(h)=\{v\}$, $v$ is not a range of $h$ and $v$ is not a source or range of a hyperedge $g\neq h$.
\end{deff}

\begin{deff}
Suppose that $H$ is basic and let $v\in H^0$ be a lonely vertex. Let $H'$ be the basic Bergman graph obtained from $H$ by removing $v$ and the unique hyperedge $h$ having $v$ as a source. Then we say that $H'$ can be obtained from $H$ by {\it lonely vertex elimination}.
\end{deff}

\begin{thm}\label{thmbaseelim}
Suppose that $H$ is basic and that the Bergman graph $H'$ can be obtained from $H$ by lonely vertex elimination. Then $B(H')\mor B(H)$.
\end{thm}
\begin{proof}
By Remark \ref{rembergalg}(d), $B(H)\cong L(H)$ and $B(H')\cong L(H')$. It follows from Definition \ref{defhlpa} that there is a $K$-algebra homomorphism $\phi:L(H')\to L(H)$ such that $\phi(u)=u$, $\phi(g_{ij})=g_{ij}$ and $\phi(g_{ij}^*)=g_{ij}^*$ for any $u\in (H')^0$, $g\in (H')^1$, $i\in I_g$ and $j\in J_g$. We will show that $\phi$ is injective, that $\phi(L(H'))=eL(H)e$ where $e=\sum_{u\in H^0\setminus\{v\}}u$, and that the idempotent $e$ is full in $L(H)$ (i.e. $L(H)=L(H)eL(H)$). The assertion of the theorem will then follow.

That $\phi$ is injective follows from Theorem \ref{thmbasis} since $\phi$ maps distinct basis elements to distinct basis elements. Next we show that $\phi(L(H'))=eL(H)e$. Clearly $\phi(L(H'))\subseteq eL(H)e$. Let now $a\in eL(H)e$. Then $a$ can be written as a $K$-linear combination of basis paths that neither start nor end in $v$. Let $p$ be one of these basis paths. Clearly $v$ cannot be a letter of $p$. Suppose $h_{ij}$ is a letter of $p$, where $i\in I_h$ and $j\in J_h$. Since $I_h$ has only one element (namely $(v,1)$), $i=i_h$. Moreover the source of the edge $h_{i_hj}$ in the graph $\widehat E$ is $v$. Hence $h_{i_hj}$ cannot be the first letter of $p$. But the only edges in $\widehat E$ ending in $v$ are the edges $h_{i_hj'}^*$ where $j'\in J_h$. Hence $p$ contains a subword $h_{i_hj'}^*h_{i_hj}$ where $j,j'\in J_h$. But this contradicts the assumption that $p$ is a basis path. Similarly one can show that none of the letters of $p$ equals $h_{ij}^*$, where $i\in I_h$ and $j\in J_h$. It follows that $p\in\phi(L(H'))$. Thus we have shown that $\phi(L(H'))=eL(H)e$.

It remains to show that $e$ is full. Clearly the ideal $L(H)eL(H)$ contains all generators of $L(H)$ except $v$ (since all the edges in $\widehat E$ start or end in a vertex not equal to $v$). But $v=\sum_{j\in J_h}h_{i_hj}h^*_{i_hj}$ by relation (iii) in Definition \ref{defhlpa}. Hence $L(H)eL(H)$ also contains $v$ and therefore $L(H)=L(H)eL(H)$.

We have shown that $L(H')$ is isomorphic to the corner $eL(H)e$ of $L(H)$, and that $e$ is a full idempotent. It follows that $L(H')\mor L(H)$, see for example the introduction of \cite{abrams-ruiz-tomforde}.
\end{proof}

\begin{deff}\label{deflonep}
We call a generator $x\in X$ {\it lonely} if there is an $i\in I$ such that $x=a_i$, $x$ is not a summand of $b_i$ and $x$ is not a summand of $a_j$ or $b_j$ for any $j\in I\setminus\{i\}$. 
\end{deff}

\begin{deff}
Suppose that $(X,R)$ is basic and let $x\in X$ be a lonely generator. Let $(X',R')$ be the Bergman presentation obtained from $(X,R)$ by removing the generator $x$ and the unique relation $a_i=b_i$ where $x=a_i$. Then we say that $(X',R')$ can be obtained from $(X,R)$ by {\it lonely generator elimination}.
\end{deff}

\begin{thm}\label{thmbaseelimp}
Suppose that $(X,R)$ is basic and that the Bergman presentation $(X',R')$ can be obtained from $(X,R)$ by lonely generator elimination. Then $B(X',R')\mor B(X,R)$.
\end{thm}
\begin{proof}
The theorem follows from Theorem \ref{thmbaseelim}.
\end{proof}

\begin{example}
Let $H$ and $H'$ be the Bergman graphs
\[H:\,\vcenter{\vbox{\xymatrix{
&f\ar@[red]@{->}[dr]\ar@[red]@{->}[dl]&\\
\mathcircled{u}\ar@/^1.5pc/@[red]@{-}[ur]&g\ar@[red]@{-}[r]\ar@[red]@{->}[l]&\mathcircled{w}\\
&h\ar@[red]@{->}[ur]\ar@/_1.5pc/@[red]@{->}[ur]\ar@[red]@{->}[ul]&\\
&\mathcircled{v}\ar@[red]@{-}[u]&
}}}
\quad\quad \text{ and }\quad\quad H':\,\vcenter{\vbox{\xymatrix{
&f\ar@[red]@{->}[dr]\ar@[red]@{->}[dl]&\\
\mathcircled{u}\ar@/^1.5pc/@[red]@{-}[ur]&g\ar@[red]@{-}[r]\ar@[red]@{->}[l]&\mathcircled{w}
}}}.\]
Then $B(H)\mor B(H)$ by Theorem \ref{thmbaseelim}.
\end{example}

\subsection{Collapsing}

\begin{deff}\label{defcollapsp}
Suppose that $(X,R)$ is basic. Let $x\in X$ and $i\in I$ such that $x=a_i$ and $x$ is not a summand of $b_i$. For any $j\in I\setminus\{i\}$ let $a'_j$ (resp. $b'_j$) be the element of $\l X\r $ obtained from $a_{j}$ (resp. $b_{j}$) by replacing in the expression of $a_j$ (resp. $b_j$) as a sum of generators any occurrence of $x$ by $b_{i}$. Let $(X',R')$ be the basic Bergman presentation obtained from $(X,R)$ by removing the $i$-th relation $a_{i}=b_{i}$ and the generator $x$, and replacing each of the other relations $a_j=b_j$, where $j\in I\setminus\{i\}$, by the corresponding relation $a'_j=b'_j$. Then we say that $(X',R')$ can be obtained from $(X,R)$ by {\it collapsing}.
\end{deff}

\begin{lemma}\label{lemcollapsp}
Suppose that $(X,R)$ is basic and that the Bergman presentation $(X',R')$ can be obtained from $(X,R)$ by collapsing. Then $(X',R')$ can be obtained from $(X,R)$ by a finite number of red shift moves followed by a lonely generator elimination.
\end{lemma}
\begin{proof}
We keep the notation from Definition \ref{defcollapsp}. Let $(X,R'')$ be the Bergman presentation obtained from $(X,R)$ by replacing each of the relations $a_j=b_j$, where $j\in I\setminus\{i\}$, by the corresponding relation $a'_j=b'_j$. Then $(X,R'')$ can be obtained from $(X,R)$ by a finite number of red shift moves. Clearly $x$ is lonely with respect to the Bergman presentation $(X,R'')$, and $(X',R')$ can be obtained from $(X,R'')$ by a lonely generator elimination. 
\end{proof}

\begin{thm}\label{thmcollapsp}
Suppose that $(X,R)$ is basic and that the Bergman presentation $(X',R')$ can be obtained from $(X,R)$ by collapsing. Then $B(X',R')\mor B(X,R)$.
\end{thm}
\begin{proof}
The theorem follows from Lemma \ref{lemcollapsp}, Theorem \ref{thmshiftredp} and Theorem \ref{thmbaseelimp}.
\end{proof}


\begin{deff}
Suppose that $H$ is basic. Let $v\in H^0$ and $h\in H^1$ such that $s(h)=\{v\}$ and $v$ is not a range of $h$. Let $H'$ be the the basic Bergman graph obtained from $H$ by removing the hyperedge $h$ and the vertex $v$, and replacing $s$ and $r$ by $s'$ and $r'$, respectively, which are defined by
\begin{align*}
s'(g)(u)=\begin{cases}
s(g)(u),\quad&\text{if } u\not\in r(h),\\
s(g)(u)+s(g)(v),\quad&\text{if }u\in r(h),
\end{cases}\quad\text{ and }\quad r'(g)(u)=\begin{cases}
r(g)(u),\quad&\text{if } u\not\in r(h),\\
r(g)(u)+r(g)(v),\quad&\text{if }u\in r(h).
\end{cases}
\end{align*}
Then we say that $H'$ can be obtained from $H$ by {\it collapsing}.
\end{deff}

\begin{lemma}\label{lemcollaps}
Suppose that $H$ is basic and that the Bergman graph $H'$ can be obtained from $H$ by collapsing. Then $H'$ can be obtained from $H$ by a finite number of red shift moves followed by a lonely vertex elimination.
\end{lemma}
\begin{proof}
See the proof of Lemma \ref{lemcollapsp}.
\end{proof}

\begin{thm}\label{thmcollaps}
Suppose that $H$ is basic and that the Bergman graph $H'$ can be obtained from $H$ by collapsing. Then $B(H')\mor B(H)$.
\end{thm}
\begin{proof}
The theorem follows from Theorem \ref{thmcollapsp} (or alternatively from Lemma \ref{lemcollaps}, Theorem \ref{thmshiftred} and Theorem \ref{thmbaseelim}).
\end{proof}

\begin{example}
Let $H$ and $H'$ be the Bergman graphs
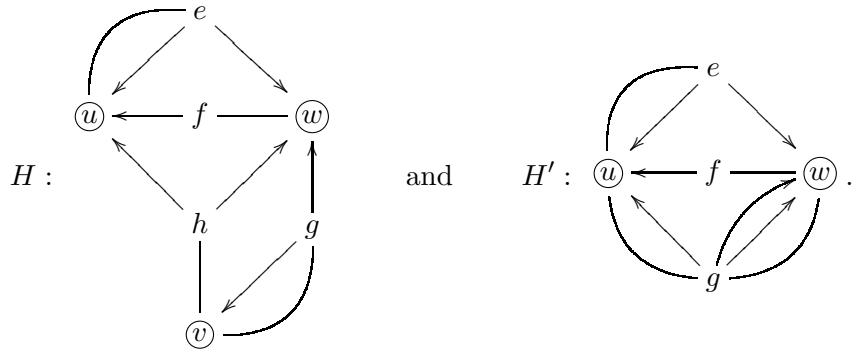
\[H:\,\vcenter{\vbox{\xymatrix{
&e\ar@[red]@{->}[dr]\ar@[red]@{->}[dl]&\\
\mathcircled{u}\ar@/^1.5pc/@[red]@{-}[ur]&f\ar@[red]@{-}[r]\ar@[red]@{->}[l]&\mathcircled{w}\\
&h\ar@[red]@{->}[ur]\ar@[red]@{->}[ul]&g\ar@[red]@{->}[u]\ar@[red]@{->}[dl]\ar@/^1.5pc/@[red]@{-}[dl]\\
&\mathcircled{v}\ar@[red]@{-}[u]&
}}}
\quad\quad \text{ and }\quad\quad H':\,\vcenter{\vbox{\xymatrix{
&e\ar@[red]@{->}[dr]\ar@[red]@{->}[dl]&\\
\mathcircled{u}\ar@/^1.5pc/@[red]@{-}[ur]&f\ar@[red]@{-}[r]\ar@[red]@{->}[l]&\mathcircled{w}\\
&g\ar@[red]@{->}[ur]\ar@/^0.8pc/@[red]@{->}[ur]\ar@[red]@{->}[ul]\ar@/_1.2pc/@[red]@{-}[ur]\ar@/^1.2pc/@[red]@{-}[ul]&
}}}.\]
Then $B(H)\mor B(H')$ by Theorem \ref{thmcollaps}. Note that $H'$ can be obtained from $H$ by a red shift move followed by a lonely vertex elimination.
\end{example}

\subsection{Insplitting}

\begin{deff}\label{definsplitp}
Suppose that $(X,R)$ is basic. Let $x_1\in X$ and $i\in I$ such that $x_1=a_i$ and there is no  $j\in I\setminus\{i\}$ such that $x_1$ is a summand of $a_j$. For any $j\in I$ let $n_j$ be the number of occurrences of $x_1$ in the expression of $b_j$ as a sum of generators. Suppose that the set $S=\{(j,k)\mid j\in I,~1\leq k\leq n_j\}$ is not the empty set and let $S=S_1\sqcup \dots\sqcup S_t$ be a partition of $S$ into pairwise disjoint nonempty subsets. Let $\phi:S\to \{1,\dots,t\}$ be the map defined by $\phi(j,k)=p$ for any $(j,k)\in S_p$ where $1\leq p\leq t$. Let $(X',R')$ be the basic Bergman presentation obtained from $(X,R)$ by
\begin{itemize}
\item adding generators $x_{2},\dots, x_t$,
\item replacing for any $j\in I\setminus\{i\}$ the $j$-th relation $a_{j}=b_{j}$ by the relation $a_j=b'_j$, where $b'_j$ is obtained from $b_j$ by replacing the summand $n_jx_1$ by the sum $ x_{\phi(j,1)}+\dots+x_{\phi(j,n_j)}$, and 
\item replacing the $i$-th relation $x_1=b_i$ by the relations $x_{p}=b'_i~(1\leq p\leq t)$ where $b'_i$ is obtained from $b_i$ by replacing the summand $n_ix_1$ by the sum $x_{\phi(i,1)}+\dots+x_{\phi(i,n_i)}$.
\end{itemize} 
Then we say that $(X',R')$ can be obtained from $(X,R)$ by {\it insplitting}.
\end{deff}

\begin{lemma}\label{leminsplitp}
Suppose that $(X,R)$ is basic and that the Bergman presentation $(X',R')$ can be obtained from $(X,R)$ by insplitting. Then $(X,R)$ can be obtained from $(X', R')$ by $t-1$ red shift moves followed by $t-1$ collapsings.
\end{lemma}
\begin{proof}
We keep the notation from Definition \ref{definsplitp}. Let $(X',R'')$ be the Bergman presentation obtained from $(X',R')$ by replacing for any $2\leq p\leq t$ the relation $x_{p}=b'_i$ by the relation $x_{p}=x_{1}$. Clearly $(X',R'')$ can be obtained from $(X',R')$ by $t-1$ red shift moves. Let $(X'',R''')$ be the Bergman presentation obtained from $(X',R'')$ by removing the generators $x_{p}$ and the relations $x_{p}=x_{1}$ where $2\leq p\leq t$, and by replacing any occurrence of a generator $x_{p}~(2\leq p\leq t)$ in the remaining relations by $x_{1}$. Clearly $(X'',R''')$ can be obtained from $(X',R'')$ by $t-1$ collapsings. One checks easily that $X''=X$ and $R'''=R$.
\end{proof}

\begin{thm}\label{thminsplitp}
Suppose that $(X,R)$ is basic and that the Bergman presentation $(X',R')$ can be obtained from $(X,R)$ by insplitting. Then $B(X',R')\mor B(X,R)$.
\end{thm}
\begin{proof}
The theorem follows from Lemma \ref{leminsplitp}, Theorem \ref{thmshiftredp} and Theorem \ref{thmcollapsp}.
\end{proof}


\begin{deff}
Suppose that $H$ is basic. Let $v_1\in H^0$ and $h_1\in H^1$ such that $s(h_1)=\{v_1\}$ and there is no $g\in H^1\setminus\{h_1\}$ such that $v_1$ is a source of $g_1$. Suppose that the set $S=\{(g,k)\mid g\in H^1,1\leq k\leq r(g)(v_1)\}$ is not the empty set and let $S=S_1\sqcup \dots\sqcup S_t$ be a partition of $S$ into pairwise disjoint nonempty subsets. Let $\phi:S\to \{1,\dots,t\}$ be the map defined by $\phi(g,k)=p$ for any $(g,k)\in S_p$ where $1\leq p\leq t$. Let $H'$ be the basic Bergman graph obtained from $H$ by adding vertices $v_2,\dots, v_t$, hyperedges $h_{2},\dots,h_t$, and by replacing $s$ and $r$ by $s'$ and $r'$, respectively, which are defined by
\begin{align*}
&s'(g)=s(g),\quad r'(g)(u)=r(g)(u),\quad r'(g)(v_{q})=\#\{(g,k)\mid 1\leq k\leq r(g)(v_1),~\phi(g,k)=q\},\\
&s'(h_{p})=v_{p},\quad r'(h_{p})(u)=r(h)(u),\quad  r'(h_{p})(v_{q})=
\#\{(h,k)\mid 1\leq k\leq r(h)(v_1),~\phi(h,k)=q\}
\end{align*}
for any $g\in H^1\setminus\{h_1\}$, $u\in H^0\setminus\{v_1\}$ and $1\leq p,q\leq t$.
Then we say that $H'$ can be obtained from $H$ by {\it insplitting}.
\end{deff}

\begin{lemma}\label{leminsplit}
Suppose that $H$ is basic and that the Bergman graph $H'$ can be obtained from $H$ by insplitting. Then $H$ can be obtained from $H'$ by $t-1$ red shift moves followed by $t-1$ collapsings.
\end{lemma}
\begin{proof}
See the proof of Lemma \ref{leminsplitp}.
\end{proof}

\begin{thm}\label{thminsplit}
Suppose that $H$ is basic and that the Bergman graph $H'$ can be obtained from $H$ by insplitting. Then $B(H')\mor B(H)$.
\end{thm}
\begin{proof}
The theorem follows from Theorem \ref{thminsplitp} (or alternatively from Lemma \ref{leminsplit}, Theorem \ref{thmshiftred} and Theorem \ref{thmcollaps}).
\end{proof}

\begin{example}
Let $H$ and $H'$ be the Bergman graphs
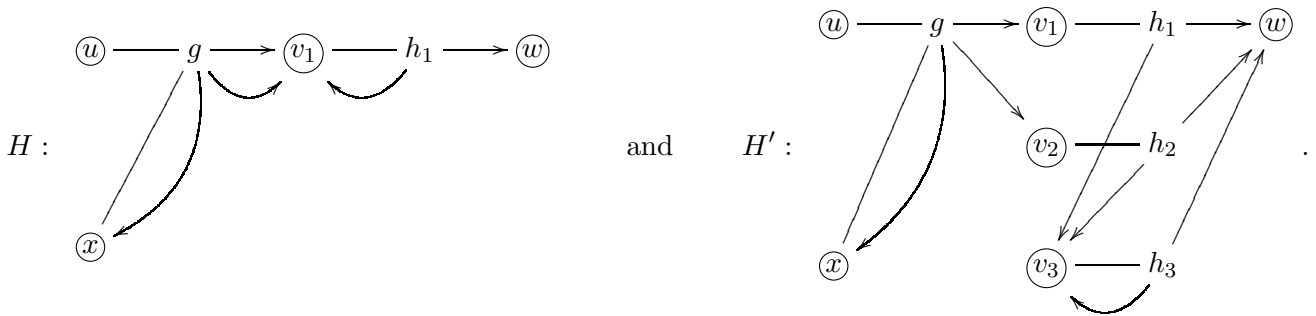
\[H:\,\,\vcenter{\vbox{\xymatrix{
\mathcircled{u}\ar@[red]@{-}[r]&g\ar@[red]@{-}[ddl]\ar@/^1.5pc/@[red][ddl]\ar@[red][r]\ar@/_1.5pc/@[red][r]&\mathcircled{v_1}\ar@[red]@{-}[r]&h_1\ar@[red]@{->}[r]\ar@[red]@/^1.5pc/@{->}[l]&\mathcircled{w}\\
&&&&\\
\mathcircled{x}&&&&
}}}
\quad\quad \text{ and }\quad\quad H':\,\,\vcenter{\vbox{\xymatrix{
\mathcircled{u}\ar@[red]@{-}[r]&g\ar@[red]@{-}[ddl]\ar@/^1.5pc/@[red][ddl]\ar@[red][r]\ar@[red][dr]&\mathcircled{v_1}\ar@[red]@{-}[r]&h_1\ar@[red]@{->}[r]\ar@[red]@{->}[ddl]&\mathcircled{w}\\
&&\mathcircled{v_2}\ar@[red]@{-}[r]&h_2\ar@[red]@{->}[ur]\ar@[red]@{->}[dl]&\\
\mathcircled{x}&&\mathcircled{v_3}\ar@[red]@{-}[r]&h_3\ar@[red]@{->}[uur]\ar@[red]@/^1.5pc/@{->}[l]&
}}}.\]
$~$\\
\\
Then $B(H)\mor B(H')$ by Theorem \ref{thminsplit} (choose the partition $S=S_1\sqcup S_2\sqcup S_3$ where $S_1=\{(g,1)\}$, $S_2=\{(g,2)\}$ and $S_3=\{(h,1)\}$). Note that $H$ can be obtained from $H'$ by two red shift moves followed by two collapsings.
\end{example}

\section{Connections to Tietze transformations}

Let $(X,R)$, where $R=\{a_i=b_i\mid~i\in I\}$, be an abelian monoid presentation. 
We call a relation $a_i=b_i$ {\it superfluous} if the equality $a_i=b_i$ holds in the abelian monoid $\l X\mid R'\r$ where $R'=\{a_j=b_j\mid~j\in I\setminus\{i\}\}$. 

\begin{deff}\label{deftietze}
Let $(X,R)$ and $(X',R')$ be good abelian monoid presentations. Suppose that one of the conditions (A)-(D) below is satisfied.
\begin{enumerate}[(A)]
\item $X'$ is obtained from $X$ by adding a generator $x\not\in X$, and $R'$ is obtained from $R$ by adding a relation $x=b$ where $b\in \l X\r$. 
\medskip
\item $X'$ is obtained from $X$ by removing a generator $x\in X$, and $R'$ is obtained from $R$ by removing a relation $x=b$ where $b\in \l X\setminus\{x\}\r$.
\medskip
\item $X'=X$ and $R'$ is obtained from $R$ by adding a superfluous relation.
\medskip
\item $X'=X$ and $R'$ is obtained from $R$ by removing a superfluous relation.
\end{enumerate}
Then we say that $(X',R')$ can be obtained from $(X,R)$ by a {\it Tietze transformation}.
\end{deff}

\begin{thm}\label{thmtietze}
Let $(X,R)$ and $(X',R')$ be finite and good abelian monoid presentations. Then $\langle X\mid R\rangle\cong\langle X'\mid R'\rangle$ if and only if $(X',R')$ can be obtained from $(X,R)$ by a finite sequence of Tietze transformations.
\end{thm}
\begin{proof}
The proof is essentially the same as the proof of \cite[Proposition 2.1 in Chapter II ]{lyndon-schupp}.
\end{proof}


\begin{cor}\label{cortietze2}
Each of the moves in Section 4 is a composition of a finite number of Tietze transformations (if one neglects the colouring of the relations in Bergman presentations).
\end{cor}
\begin{proof}
Let $(X,R)$ and $(X',R')$ be Bergman presentations and suppose that $(X',R')$ can be obtained from $(X,R)$ by one of the moves defined in Section 4. Then $B(X,R)\mor B(X',R')$ as shown in Section 4. It follows from Lemma \ref{lemvmon} that $\l X\mid R\r\cong\V(B(X,R))\cong \V(B(X',R'))\cong\l X'\mid R'\r$. Thus, by Theorem \ref{thmtietze}, $(X',R')$ can be obtained from $(X,R)$ by a finite sequence of Tietze transformations.
\end{proof}

\begin{example}
Let $(X,R)$ and $(X',R')$ be Bergman presentations and suppose that $(X',R')$ can be obtained from $(X,R)$ by a red shift move. Then there is an $i\in I$ and $a,b\in  \l X\r$ such that $a_{i}=a$ and $b_i=b$ in the abelian monoid $\l X\mid a_{j}=b_{j}~(j\in I\setminus{i})\r$, and $(X',R')$ is obtained from $(X,R)$ by replacing the $i$-th relation $a_{i}=b_{i}$ by the relation $a=b$. One can factor this move into two Tietze transformations: First add the superfluous relations $a=b$, then remove the superfluous relation $a_i=b_i$. Hence any red shift move is a composition of a Tietze transformation of type (C) and one of type (D).
\end{example}

\begin{rmk}
For basic Bergman presentations (which are essentially just the finite and good abelian monoid presentations), the Tietze transformations of type (B) in Definition \ref{deftietze} are precisely the lonely generator eliminations. Hence, by Theorem \ref{thmbaseelimp}, if a basic Bergman presentation $(X',R')$ can be obtained from another one $(X,R)$ by a Tietze transformation of type (A) or (B), then $B(X,R)\mor B(X',R')$. If the Tietze transformations of type (C) and (D) had the same property (i.e. if they preserved the Morita equivalence class on the level of algebras), then it would follow from Lemma \ref{lemvmon} and Theorem \ref{thmtietze} that $B(X,R)\mor B(X',R')\Leftrightarrow \l X\mid R\r\cong\l X'\mid R'\r$ for any basic Bergman presentations $(X,R)$ and $(X',R')$. Unfortunately, the Tietze transformations of type (C) and (D) do in general not preserve the Morita equivalence class on the level of algebras. For example, the basic Bergman presentation $(\{x\},\{x=x\})$ can be obtained from $(\{x\},\emptyset)$ by a Tietze transformation of type (C), but $B(\{x\},\{x=x\})\cong K[x,x^{-1}]$ and $B(\{x\},\emptyset)\cong K$ are not Morita equivalent.
\end{rmk}


\begin{thebibliography}{9999}

\bibitem{lpabook} G. Abrams, P. Ara, M. Siles Molina, Leavitt path algebras. Lecture Notes in Mathematics, vol. 2191, Springer Verlag, 2017.

\bibitem {AA} G. Abrams, G. Aranda Pino, \emph{The Leavitt path algebra of a graph}, J. Algebra {\bf 293} (2005) 319--334.

\bibitem {AN} G. Abrams, T.G. Nam, \emph{Corners of Leavitt path algebras of finite graphs are Leavitt path algebras}, J. Algebra {\bf 547} (2020), 494--518.

\bibitem{abrams-ruiz-tomforde} G. Abrams, E. Ruiz, M. Tomforde, \emph{Morita equivalence for graded rings}, J. Algebra {\bf 617} (2023), 79--112.



\bibitem{aragoodearl} P. Ara, K. Goodearl, \emph{Leavitt path algebras of separated graphs}, J. reine angew. Math. {\bf 669} (2012), 165--224.


\bibitem {AMP} P. Ara, M.A. Moreno, E. Pardo, \emph{Nonstable K-thory for Graph Algebras}, Algebr. Represent. Theory {\bf 10} (2007) 157--168.



\bibitem{bergman74} G.M. Bergman, \emph{Coproducts and some universal ring constructions}, Trans. Amer. Math. Soc. {\bf 200} (1974), 33--88.










\bibitem {H-1}R. Hazrat, \emph{The graded structure of Leavitt path algebras}, Israel
J. Math,  {\bf 195} (2013) 833--895.



\bibitem{kbgm} M. Kanuni, D. Martin Barquero, C. Martin Gonzalez, M. Siles Molina, \emph{Classification of Leavitt path algebras with two vertices}, Mosc. Math. J. {\bf 19} (2019), no. 3, 523--548.

\bibitem{lam} T.Y. Lam, Lectures on modules and rings. Graduate Texts in Mathematics, vol. 189, Springer Verlag, 1999.

\bibitem {vitt56} W.G. Leavitt, \emph{Modules over rings of words}, Proc. Amer. Math. Soc. {\bf 7} (1956) 188--193. 

\bibitem {vitt57} W.G. Leavitt, \emph{Modules without invariant basis number,} Proc. Amer. Math. Soc. {\bf 8} (1957) 322--328.

\bibitem {vitt62} W.G.  Leavitt, \emph{The module type of a ring}, Trans. Amer. Math. Soc. {\bf 103} (1962) 113--130. 

\bibitem {vitt65} W.G. Leavitt, \emph{The module type of homomorphic images,} Duke Math. J. {\bf 32} (1965) 305--311.

\bibitem{lyndon-schupp} R.C. Lyndon, P.E. Schupp, Combinatorial group theory. Springer Verlag, 1977.

\bibitem{mohan-suhas} R. Mohan, B.N. Suhas, \emph{Cohn-Leavitt path algebras of bi-separated graphs}, Comm. Algebra {\bf 49} (2021), 1991--2021. 




\bibitem{preusser-1} R. Preusser, \emph{The V-monoid of a weighted Leavitt path algebra}, Israel J. Math. {\bf 234} (2019), no. 1, 125--147.

\bibitem{Raimund2} R. Preusser, \emph{Leavitt path algebras of hypergraphs},  Bull. Braz. Math. Soc. (N.S.) {\bf 51} (2020), 185--221. 





\end{thebibliography}
\end{document}